\documentclass{amsart}

        \usepackage {amssymb,amsthm,amsxtra,amscd,amsmath,amsfonts}

        \newcommand{\K}{\ensuremath{\Bbbk}}
        \newcommand{\N}{\ensuremath{\mathbb{N}}}
        \newcommand{\Q}{\ensuremath{\mathbb{Q}}}
        \newcommand{\Z}{\ensuremath{\mathbb{Z}}}

        \newcommand{\degr}[1]{\textup{deg}(#1)}
        \newcommand{\diff}{\ensuremath{\mathit{diff}}}

        \newcommand{\Aut}{\textup{Aut}}
        \newcommand{\Hom}{\textup{Hom}}
        \newcommand{\grHom}{\textup{grHom}}

        \newcommand{\del}{\ensuremath{\partial}}
        \newcommand{\onto}{\twoheadrightarrow}

	\newcommand{\Renvelope}{\ensuremath{R \otimes R^o}}
	\newcommand{\Dnaught}{\ensuremath{D^0_{\K}(R)}}

        \theoremstyle{plain}
                \newtheorem{theorem}{Theorem}[section]
                \newtheorem{lemma}{Lemma}[section]
                \newtheorem{corollary}{Corollary}[section]
                \newtheorem{proposition}{Proposition}[section]
                \newtheorem{remark}{Remark}[section]
                \newtheorem{definition}{Definition}[section]
                \newtheorem*{theorem*}{Theorem}
                \newtheorem*{remark*}{Remark}
        \theoremstyle{definition}
                
        \numberwithin{equation}{subsection}
        
        \newcommand{\ignore}[1]{}
        \newcommand{\term}[1]{\textit{#1}}
        \newcommand{\mynote}[1]{}

\begin{document}

        \title[Differential Operators on the free algebras]
	{Differential Operators on the free algebras}
        \author[Iyer]{Uma N. Iyer} 
        \author[McCune]{Timothy C. McCune} 
        \email{uma.iyer@bcc.cuny.edu, tcmccune@yahoo.com}
        \address{Department of Mathematics and Computer Science,
		Bronx Community College,
		2155 University Avenue, Bronx, New York 10453.
		}

\begin{abstract}
Following the definitions of the algebras of differential operators,
$\beta$-differential operators, and the quantum differential operators
on a noncommutative (graded) algebra given in \cite{LR}, we describe
these operators on the free associative algebra. We further study
their properties.

\vspace{0.2in}
\end{abstract}
\maketitle              
\section{Introduction}

Let $\K$ denote a field, and $R$ a $\K$-algebra.  Any $R$-bimodule $M$
under consideration
satisfies 
$\alpha m = m \alpha$ for $m\in M$ and $\alpha \in \K$. 
In \cite{LR} V.A.Lunts and
A.L.Rosenberg defined an $R$-subbimodule,
\term{the differential part of} $M$, denoted by $M_{\diff}$.
The $R$-bimodule $M_{\diff}$ has a filtration $M_0 \subset M_1 \subset \cdots$,
where $M_i$ is called \term{the $i$-th differential part of} $M$.
In particular, the elements of the
differential part of the $R$-bimodule $\Hom _{\K}(R,R)$
are called \term{the differential operators on} $R$. 
This $R$-bimodule is an algebra and is denoted by $D_{\K}(R)$.  
The differential operators on $R$ of \term{order} $\leq n$ are elements
of the $n$-th differential part of $\Hom _{\K}(R,R)$. Analogously, 
$\beta$-differential operators for any $\Gamma$-graded $\K$ algebras
($\Gamma$ is an abelian group) are defined in \cite{LR}. 
In the particular case of $\Gamma = \mathbb{Z}_2$, the $\beta$-differential
operators are the super differential operators. In general, under
certain conditions on $\beta$, the $\beta$-differential operators
are the coloured differential operators on coloured algebras.
A further notion of quantum differential operators on a $\Gamma$-graded
$\K$ algebra is defined, which allows for viewing action of a quantum group
on a ring through its Hopf structure to be via quantum differential operators.

The differential operators on the polynomial algebra 
and the supercommutative free algebra over fields of characteristic $0$
have been widely studied. Smith in \cite{S} 
studied differential operators on the 
affine and projective lines for nonzero characteristic.
In particular, the algebra of
differential operators in the case of nonzero characteristic is not finitely
generated. 

In this  paper we investigate the algebra of differential operators on
the free algebra and the $\beta$ and quantum 
differential operators on the free 
algebra graded by $\Z^n$. 

The preliminaries are given in the section \ref{S:Prelim}.   
The algebra of usual differential operators on the 
free algebras are described and studied in section
\ref{S:free}.  When $R$ is the free algebra over $n$ variables,
$x_1,\cdots ,x_n$, the algebra of
$0$-th order differential operators, $D_{\K}^0 (R)$, is generated by
left and right multiplication homomorphisms. 
The $D_{\K}^0(R)$-module of the first order differential operators 
is generated by derivations. Higher order differential operators
are defined in the definition \ref{Def:usualdef}. Unlike the case
of polynomial algebras, the first order differential operators
do not generate the algebra of differential operators, $D_{\K}(R)$,
even when the characteristic of the underlying field is $0$. We are able
to prove that $D_{\K}(R)$ is a simple algebra. We also describe several
properties of these new operators in this section. Of particular interest
is the fact that every differential operator can be written as a power
series in a unique way (Remark \ref{Remark:canonical-series}).

An analogous study of the $\beta$-differential operators is
pursued in section \ref{S:betafree}. Since the $\beta$-differential
operators are similar to the usual differential operators, the proofs of
the results are not presented.

Finally, a study of the 
quantum differential operators on the free algebras
is presented in section \ref{S:qfree}.
We present some properties of these quantum differential operators 
in this section. The algebra of quantum differential operators
on the free algebra, $D_q(R)$, has a more complicated structure than 
$D_{\K}(R)$. A study of the algebra structure of $D_q(R)$ (for
example, checking whether it is simple, or a domain) 
will be pursued in a later project. In this paper, we primarily 
address the properties of $D_{\K}(R)$.

There are no
restrictions on the characteristic of the underlying field 
unless otherwise mentioned.
\subsection*{Acknowledgements}
The first author would like to thank Professor S.P.Smith
for his help during the preparation of this work. 
She also gratefully acknowledges support from
the PSC-CUNY Grants, award \# 62279-00 40 (2009-2010) and 
award \# 63032-00 41 (2010-2011).

\section{Preliminaries}\label{S:Prelim}
Let $\K$ be a field,
and $R$ be an associative, unital $\K$-algebra.
We let $\bigotimes$ denote $\bigotimes _{\K}$.
Let $R^e$ denote
$R\bigotimes R^o$ where $R^o$ is the opposite ring of $R$.
Let $M$ be an $R$-bimodule, equivalently, a left $R^e$-module.
We recall the definition of the \term{differential part} of $M$, denoted by
$M_{\mathit{diff}}$  
from \cite{LR}.
The \term{centre} of $M$ is the $\K$-vector space
\[Z(M) = \{ m\in M | rm =mr \textit{  for all } r\in R \}.
\]
If $R$ is commutative, then $Z(M)$ is an $R^e$-submodule of $M$.
The \term{differential part}, of $M$ is
the $R$-bimodule 
$M_{\mathit{diff}} = \cup _{i=0}^{\infty} M_i$, where
$M_i$ is the $R$-bimodule generated by $\{ m \in M| \bar{m} 
\in Z(M/M_{i-1}) \}$ with $M_{-1} = 0$.  Each $M_i$ is called the
\term{$i$-th differential part }of $M$.\\
\textbf{Notation :} For $r\in R$ and $m\in M$, let
$[m,r]:= mr - rm$ and $[r,m] := rm - mr$. For any $\K$-algebra $A$, we set
$[a,b] := ab-ba$ for $a,b\in A$.

The vector space $\mathit{Hom}_{\K}(R,R)$ is an $R$-bimodule as follows:
For $r,s \in R$ and $\varphi \in \mathit{Hom}_{\K}(R,R)$, we let
$
(r \varphi) (s) = r\varphi (s)$, and $(\varphi r)(s)= \varphi (rs)$.

The $R$-bimodule $D^m_{\K}(R)$ of 
differential operators on $R$ of order $m$
is the $m$-th differential part of $\mathit{Hom}_{\K}(R,R)$.
Note that, if $R$ is commutative,
a homomorphism
$\varphi \in \mathit{Hom}_{\K}(R,R)$ is said to be a differential operator of
order $m$ if 
$[\cdots [[\varphi ,r_0],r_1],\cdots ,r_m] = 0$ for
all $r_i \in R$.

We see that $\varphi \in Z(\mathit{Hom}_{\K}(R,R))$ 
implies, $\varphi (r) = (\varphi r)(1)
= (r\varphi )(1) = r\varphi (1)$.  
For any $r\in R$, let $\lambda_r, \rho_r \in \Hom (R,R)$ be
homomorphisms given by $\lambda_r (s) = rs$ and $\rho_r (s) =sr$
for $r\in R$.
That is,, $\varphi = \rho _{\varphi (1)}$, the
homomorphism given by right multiplication by $\varphi (1)$.  
Thus, 
$
D^0_{\K}(R)$ is the algebra generated by the
set $\{ \lambda _r , \rho _r | r \in R \}$.
There is a surjection $R \otimes_{Z(R)} R^o \onto D^0_{\K} (R)$
given by
$a\otimes b^o \mapsto \lambda_a \rho_b$
for $a,b \in R$
where $Z(R)$ is the centre of $R$.

By Remark 1.1.2.8 of \cite{LR}, we have 
$D^i_{\K}(R) D^j_{\K}(R) \subset D^{i+j}_{\K}(R)$.  Hence,
$D^0 _{\K}(R)$ is a ring, called the
ring of inner differential operators which
contains $R$ seen as $\lambda _r $ for $r\in R$.
The ring of differential operators
$D_{\K}(R)$ is the union of $D^i_{\K}(R)$.

If $\varphi$ is a derivation, then $\varphi (rs) = r\varphi (s) +
   \varphi (r) s$.  Thus, $[\varphi ,r ] = \lambda _{\varphi (r)}$. Hence,
a derivation is a first order differential operator.

For $I$ be a two sided ideal of a $\K$-algebra $R$ let 
\[
\mathcal{S}_I = \{ \varphi \in 
D_{\K}(R) \mid \varphi (I) \subset I \}\quad \textit{ and }
\quad
\mathcal{J}_I 
= \{ \varphi \in 
D_{\K}(R) \mid \varphi (R) \subset I \}.
\] 
Then $\mathcal{S}_I$ is a filtered subalgebra of $D_{\K}(R)$ 
with $\mathcal{J}_I$ as a filtered ideal. 
\begin{proposition}\label{S/RtoD}
The natural map of $\K$-algebras 
$\mathcal{S}_I/\mathcal{J}_I\to \Hom_{\K}(R/I , R/I)$ gives a 
map of filtered $\K$-algebras $\mathcal{S}_I/\mathcal{J}_I \to D_{\K}(R/I)$.
\end{proposition}
\begin{proof}
For $\varphi \in \Hom_{\K}(R,R)$, with $\varphi (I)\subset I$ 
let 
$\overline{\varphi}$ denote the corresponding map 
on $R/I$. Let $\overline{a}$ denote the image of $a\in R$.
Note, $\overline{\lambda_a}= \lambda_{\overline{a}}$,
$\overline{\rho_a}= \rho_{\overline{a}}$, and
$\overline{[\varphi ,a]} = [\overline{\varphi}, \overline{a}]$. 
The claim then follows.
\end{proof}
We now refer to \cite{LR} and present some preliminaries on
beta-differential operators which are generalizations of superderivations
and their higher powers, see \cite{L}.
Let $\Gamma$ be an abelian group. Fix a
bicharacter $\beta : \Gamma \times \Gamma \longrightarrow \K^*$.
Let $R$ be a
$\Gamma$-graded $\K$-algebra and $M$ a $\Gamma$-graded $R$-bimodule.
Let $Z_{\beta}(M)$ denote the $\beta$-{\it center} of $M$ defined as
the $\K$-span of homogeneous elements $m \in M$ such that
\[
mr = \beta (d_m , d_r) rm
\;\text{ for any homogeneous $r\in R$}, \]
where $d_x$ denotes the
degree of $x$.
The $\beta$-\emph{differential part} of $M$ is 
the $R$-bimodule 
$M_{\beta} = \cup _{i=0}^{\infty} M_{\beta ,i}$, where
$M_{\beta, i}$ is the $R$-bimodule generated by the set

\noindent
$\{ m \in M| \bar{m} 
\in Z_{\beta}(M/M_{\beta ,i-1}) \}$ 
with $M_{\beta ,-1} = 0$.  Each $M_{\beta ,i}$ is the
\term{$i$-th $\beta$-differential part }of $M$.

\noindent
\textbf{Notations:} For $r\in R$ and $m\in M$ homogeneous, let
$[m,r]_{\beta}:= mr - \beta (d_m, d_r)rm$ and 
$[r,m]_{\beta} := rm - \beta (d_r,d_m)mr$. 

For any $\Gamma$-graded
$\K$-algebra $A$, we set
$[a,b]_{\beta} := ab- \beta (d_a,d_b)ba$ for homogeneous $a,b\in A$ and
extend $[\cdot ,\cdot ]$ linearly.
Note,
\[
[ab,c]_{\beta} = \beta (d_b, d_c) [a,c]_{\beta}b + a[b,c]_{\beta} \quad 
\textit{and}\quad
[c,ab]_{\beta} = [c,a]_{\beta}b + \beta (d_c,d_a)a[c,b]_{\beta}.
\]
The $\Gamma$-graded $R$-bimodule $D^m_{\beta}(R)$ 
is the $m$-th $\beta$-differential
part of the $\Gamma$-graded $R$-bimodule $\grHom_{\K}(R,R)$.
An element of $D^m_{\beta}(R)$ is called, the $\beta$-differential operator
of order $m$.
One can see that
$
D^i_{\beta}(R) D^j_{\beta}(R) \subset D^{i+j}_{\beta}(R).
$
That is,  $D_{\beta}(R)$ is a filtered $\K$-algebra.

For a homogeneous $r\in R$, 
let $\rho^{\beta}_r \in \grHom_{\K}(R,R)$ be defined
by
\[
\rho^{\beta}_r (s) = \beta (d_r , d_s)sr.
\]
We see that $\varphi \in Z_{\beta}(\mathit{grHom}_{\K}(R,R))$ 
implies, $\varphi (r) = (\varphi r)(1)
= \beta (d_{\varphi}, d_r)r\varphi (1) $ for any homogeneous $r\in R$.
Thus, $\varphi = \rho^{\beta} _{\varphi (1)}$.
Thus, the $R$-bimodule $D^0_{\beta}(R)$ is generated by the
homomorphisms $\lambda _r, \rho^{\beta} _r$, for $r\in R$, $r$ homogeneous. 

We say that a homogenous $\varphi \in \grHom_{\K}(R,R)$ is a
left $\beta$-derivation if 
\[
\varphi (rs) = \varphi (r)s + 
	\beta (d_\varphi,d_r) r\varphi (s)
\] 
for homogeneous $r\in R$, and
any $s\in R$. Note, a left $\beta$-derivation is a 
$\beta$-differential operator
of order $1$. We see that if $\varphi$ is a homogeneous left 
$\beta$-derivation, then
$	\forall \gamma_1, \gamma_2 \in \Gamma$,
\begin{align*}
\varphi( [a,b]_{\beta}) &= [\varphi (a),b]_{\beta} + \beta (d_{\varphi},d_a) 
\left( a\varphi (b)
	-\beta (d_{\varphi}, d_a)^{-1}\beta (d_a,d_b) \varphi (b)a \right);\\
 	&= [\varphi (a),b]_{\beta} + \beta (d_{\varphi},d_a) 
        [a,\varphi (b)]_{\beta}
	\quad \textit{if } \beta (\gamma_1, \gamma_2) \beta 
        ( \gamma_2, \gamma_1) = 1 .
\end{align*}
Given two $\Gamma$-graded algebras $A,B$, the 
vector space $A\otimes B$ is a $\Gamma$-graded $\K$-algebra.  
We denote by $A\otimes ^{\beta} B$ the set $A\otimes B$, with a 
multiplication operation which reflects the $\Gamma$ grading. 
For $a,c \in A, u,v\in B$ with $b,u$ homogeneous, we have
\[
(a\otimes^{\beta} b) (u\otimes^{\beta} v) = \beta (d_b, d_u)  
(au \otimes^{\beta} bv).
\]
This multiplication also makes $A\otimes^{\beta} B$ into a $\Gamma$-graded
$\K$-algebra.
We always use this multiplication 
when we work with tensor product of $\Gamma$-graded algebras in the sections
of $\beta$-differential operators.

Given a $\Gamma$-graded algebra $R$, its $\beta$-\textit{opposite} algebra
is denoted by $R^{\beta, o}$, which as a set is $R^{\beta ,o} = 
\{ r^o \mid r\in R \}$
and the operations are $r^o+s^o = (r+s)^o$ (addition), 
$-r^o = (-r)^o$ (additive inverse), and 
$r^os^o = \beta(d_r, d_s) (sr)^o$ (multiplication for homogeneous elements).

Using these conventions we see that there is a surjection of $\Gamma$-graded
$\K$-algebras $R\otimes ^{\beta} R^{\beta ,o} \to D^0_{\beta}(R)$ given by
$a\otimes^{\beta} b^o \mapsto \lambda_a \rho_b^{\beta}$.
This surjection reduces to a surjection of $\K$-algebras
$R\otimes ^{\beta}_{Z_{\beta}(R)} R^{\beta ,o} \to D^0_{\beta}(R)$ where
$Z_{\beta}(R)$ is the $\beta$-centre of $R$; that is, it is the subalgebra
of $R$ given by
$
Z_{\beta}(R) = \{ r \in R \mid \lambda_r = \rho_r^{\beta} \}.
$

For $I$ a two sided $\Gamma$-graded ideal of $R$ let 
\[
\mathcal{S}^{\beta}_I = \{ \varphi \in 
D_{\beta}(R) \mid \varphi (I) \subset I \}\quad \textit{ and }
\quad
\mathcal{J}^{\beta}_I 
= \{ \varphi \in 
D_{\beta}(R) \mid \varphi (R) \subset I \}.
\] 
Then $\mathcal{S}^{\beta}_I$ is a $\Gamma$-graded $\Z$-filtered
subalgebra of $D_{\beta}(R)$ 
with $\mathcal{J}^{\beta}_I$ as a $\Gamma$-graded $\Z$-filtered ideal. 
The proof of the following is similar to that of Proposition \ref{S/RtoD}.
\begin{proposition}\label{super-S/RtoD}
The natural map 
$\mathcal{S}^{\beta}_I/\mathcal{J}^{\beta}_I\to \grHom_{\K}(R/I , R/I)
$ of $\Gamma$-graded $\K$-algebras 
gives a map of $\Gamma$-graded
$\Z$-filtered $\K$-algebras 
$\mathcal{S}^{\beta}_I/\mathcal{J}^{\beta}_I \to D_{\beta}(R/I)$.
\end{proposition}
Note that if $\beta \equiv 1$, then we get the usual differential operators.

Let $\Z _2:=\Z/2$. A $\Z _2$-graded $\K$-algebra
$R = R_0 \oplus R_1$ is called a \emph{superalgebra}.
Elements of $R_0$ are called \emph{even} and those of $R_1$
are called \emph{odd}.
For a homogeneous element $a\in R$,  we let $p(a)$ denote its
\emph{parity} (which is the same as \emph{degree} in this case).
For a superalgebra $R$,
we define a bicharacter $\beta :\Z_2 \times \Z_2 \longrightarrow \K^*$
by setting $\beta (x,  y) = (-1)^{xy}$ and study the $\beta$-differential
operators on a superalgebra $R$. The superdifferential operators 
(respectively, the
superderivations) are special cases of 
$\beta$-differential operators 
(respectively, the $\beta$-derivations) on
a superalgebra. More generally, when $\beta (a,b) \beta (b,a)=1$,  we get
the notion of coloured differential operators. 
Note that most of the existing studies of 
super-differential operators (respectively, coloured differential operators) 
are
on supercommutative (respectively, coloured-commutative) algebras. In this 
article
we study the $\beta$-differential operators on the free algebra on several 
variables. 

We now recall the definition of the algebra of quantum 
differential operators (\cite{LR}).
Let $\Gamma$ be an abelian group. Fix a bicharacter 
$\beta : \Gamma \times \Gamma \longrightarrow \K^*$.

Let $R$ be a $\Gamma$-graded $\K$-algebra and $M$ a $\Gamma$-graded 
$R$-bimodule.

Let $\mathcal{Z}_{q}(M)$ denote the {\it quantum-center} of $M$ defined as
	the $\K$-span of homogeneous elements $m \in M$ 
        for which there exists a $d\in \Gamma$ such that
	\[
		mr = \beta (d , d_r) rm
		\;\text{ for any homogeneous $r\in R$}.
	 \]
 For each $a\in \Gamma$, define $\sigma_a \in \grHom_{\K}(M,M)$ defined by
	$\sigma_a (m) = \beta (a,d_m)m$ for homogeneous $m\in M$, and extend 
        $\sigma_a$ linearly on all of $M$.  
	For 
	$m\in M, r\in R$, let $[m,r]_a = mr -\sigma_a(r) m$.  Using these 
        notations,
	\[
	\mathcal{Z}_{q}(M) = \K-\textit{span} \{ \textit{homogeneous } m \mid
	\exists a\in \Gamma \textit{ such that } [m,r]_a =0 \forall r\in R \}.
	\]
 	Let $M_{q,0} = R\mathcal{Z}_q(M)R$.   For $i\geq 1$, 
        $M_{q,i}$ denotes the $R$-bimodule
	generated by the set 
	\[
	 \K-\textit{span} \{ \textit{homogeneous } m \mid
	\exists a\in \Gamma \textit{ such that } [m,r]_a \in M_{q, i-1} 
                               \forall r\in R \}.
	\]
Note, $M_{q,0}\subset M_{q,1}\subset \cdots $ and 
$M_{q-\textit{diff}} = \cup_{i\geq 0}M_{q,i}$.
	When $M= \grHom_{\K}(R,R)$ we get the filtered algebra of quantum 
        differential operators 
	$D_q(R) =M_{q-\textit{diff}}$ and the $R$-bimodule of quantum 
        differential 
        operators of order $\leq i$
	is $D_{q}^i(R) = M_{q,i}$.

The algebra $D^0_q(R)$ is generated by the set 
$\{ \lambda_r, \rho_s, \sigma_a \mid r,s\in R, a\in \Gamma \}$
	where 
	\[ \lambda_r\rho_s = \rho_s\lambda_r, \hspace{0.05in}
	\sigma_a \lambda_r = \lambda_{\sigma_a(r)} \sigma_a, 
	\textit{ and }
	\sigma_a \rho_r = \rho_{\sigma_a(r)} \sigma_a.
	\]
Given a ring $S$, and a group $G$ of automorphisms of $S$, denote by $S\# G$ 
the
skew group ring on $S$ by $G$; that is, $S\# G$ is a free left $S$-module,
$\oplus_{g\in G} Sg$, with basis $G$, and with multiplication given by
\[
(r_1g_1)(r_2g_2) = r_1g_1(r_2) g_1g_2 \quad r_1,r_2\in S, g_1,g_2\in G.
\]
Now $\Gamma$ acts on $R$ via automorphisms $\sigma_a$ for $a\in \Gamma$. 
We thus get a surjection 
\[
(R\otimes_{Z(R)}R^o)\# \Gamma \to D^0_q(R) \quad \textit{given by }
(a\otimes b^o)\gamma \mapsto \lambda_a\rho_b \sigma_{\gamma}.
\]
For $i\geq 1$, each $D^i_q(R)$ is the $R$-bimodule generated by the 
$\K$-span of the set
	\[
	\{ \textit{homogeneous }\varphi \mid \exists a\in \Gamma 
        \textit{ such that } [\varphi ,r ]_a \in D^{i-1}_q(R) \}.
	\]
	Equivalently, it has been shown in \cite{IM} that each $D^i_q(R)$ is 
        the 
        $D^0_q(R)$-bimodule generated by
	the $\K$-span of the set
	\[
	\{ \textit{homogeneous }\varphi \mid  [\varphi ,r ] \in 
        D^{i-1}_q(R) \}.
	\]
	For each $a \in \Gamma$, let $\varphi \in \grHom_{\K}(R,R)$ be a left
	skew $\sigma_a$-derivation. That is, 
        $\varphi (rs) = \varphi (r)s + \sigma_a(r)\varphi(s)$ 
        $\forall r,s\in R$.
	Then $[\varphi , r]_a = \lambda_{\varphi (r)} \in D^0_q(R), 
        \forall r\in R$. 
        That is, $\varphi \in D^1_q (R)$. Similarly, for each $a\in \Gamma$, 
        let 
        $\psi \in \grHom_{\K}(R,R)$ be a right skew $\sigma_a$-derivation;
        that is,
        $\psi (rs) = \psi (r) \sigma_a(s) + r\psi (s)$ $\forall
        r,s\in R$. Then, $\forall r\in R$,
        $[\psi ,r] = \lambda_{\psi (r)}\sigma_a \in D^0_q (R)$. 
        Note that $\varphi$ is a left skew 
        $\sigma_a$-derivation, if and only if
        $\varphi \sigma_{-a}$ is a right skew $\sigma_{-a}$-derivation.

For $I$ be a two sided $\Gamma$-graded ideal of $R$ let 
\[
\mathcal{S}^{q}_I = \{ \varphi \in 
D_{q}(R) \mid \varphi (I) \subset I \}\quad \textit{ and }
\quad
\mathcal{J}^{q}_I 
= \{ \varphi \in 
D_{q}(R) \mid \varphi (R) \subset I \}.
\] 
Then $\mathcal{S}^{q}_I$ is a $\Gamma$-graded $\Z$-filtered
subalgebra of $D_{q}(R)$ 
with $\mathcal{J}^{q}_I$ as a $\Gamma$-graded $\Z$-filtered ideal. 
Again,
\begin{proposition}\label{q-S/RtoD}
The natural map 
$\mathcal{S}^{q}_I/\mathcal{J}^{q}_I\to \grHom_{\K}(R/I , R/I)$ 
of $\Gamma$-graded $\K$-algebras 
gives a map of $\Gamma$-graded
$\Z$-filtered $\K$-algebras 
$\mathcal{S}^{q}_I/\mathcal{J}^{q}_I \to D_{q}(R/I)$.
\end{proposition}
Again, if $\beta \equiv 1$, then we get the usual differential operators.

\section{The free algebra.}
\label{S:free}
Let $R= \K \langle x_1,\cdots ,x_n \rangle$ be the free
algebra over $\K$ generated by
$x_1,\cdots ,x_n$. When $n=1$, $R$ is the polynomial ring in one variable,
and the differential operators on polynomial rings have been well studied.
Therefore, {\bf assume that $n>1$.}

The ring $R$ has only associative relations.  It is a domain, and the 
monomials 
in the $x_i$ are independent over
$\K$, so we should expect to find few relations on $\Dnaught$.
For the proof, we will need the following
\begin{lemma}
  For $i = 1, \ldots, m$, let $a_i$ and $b_i$ be monic monomials.
  Let $d_i = \degr{a_i}$ and $d = \textup{max}\{ d_i \}$. 
  Then $a_i x_1^d x_2 b_i = a_j x_1^d x_2 b_j$ if and only if $a_i = a_j$
  and $b_i = b_j$.
\end{lemma}
\begin{proof}
	If $a_i = a_j$ and $b_i = b_j$, then clearly 
		$a_i x_1^d x_2 b_i = a_j x_1^d x_2 b_j$. 

	Assume $d_i \le d_j$.  Then the first $d_i$ terms of $a_i$ and $a_j$
	coincide.  Thus there is some monic monomial $c$ such that 
 	$\text{deg}(c) = d_j - d_i$ and $a_j = a_i c$.  Hence 
	  	$x_1^d x_2 b_i = c x_1^d x_2 b_j$.

	Since $\text{deg}(c) \le d$, we have $c = x_1^{d_j - d_i}$.
	Cancelling $x_1^d$ yields $x_2 b_i = x_1^{d_j - d_i} x_2 b_j$.  
	Thus $d_i = d_j$ and $a_i = a_j$.  
	Cancelling $x_2$ yields $b_i = b_j$.
\end{proof}

\begin{proposition}\label{proposition-0-order}
  The associative algebra $R\otimes R^o$ and $D^0_{\K}(R)$ are isomorphic.
\end{proposition}
\begin{proof}
	The centre of $R$ is just $\K$, and thus we have a surjection
	$R\otimes R^o \onto D^0_{\K} (R)$
	given by
	$a\otimes b^o \mapsto \lambda_a \rho_b$ for $a,b \in R$.
	It remains to show that this surjection is injective.

	An element $t$ of $\Renvelope$ may be written as 
	$t = \sum \alpha_i a_i \otimes b_i$
	where $\alpha_i \in \K$ and $a_i$ and $b_i$ are monic monomials.
	Then $t \mapsto \phi_t = \sum \alpha_i \lambda_{a_i} \rho_{b_i}$.
	Suppose $\phi_t = 0$
	Let $d_i = \degr{a_i}$ and $d = \textup{max}\{ d_i \}$.
	Then by the lemma, $\{ a_i x_1^d x_2 b_i \}$ are pairwise independent.
	Thus no two terms of $\phi_t(x_1^d x_2)$ can cancel.  
	As $R$ is the free 
	algebra, $\phi_t(x_1^d x_2)$ is 0 if and only if all $\alpha_i = 0$.
	That is, $\phi_t(x_1^d x_2) = 0$ if and only if $t = 0$.
\end{proof}

\begin{corollary}
  The Gelfand-Kirillov dimension of $D_{\K}^0(R)$ (and hence of $D_{\K}(R)$)
  is infinity.
\end{corollary}
\begin{proof}
	The Gelfand-Kirillov dimension of the free algebra $R$ is infinity.  
\end{proof}

\begin{corollary}
  The centre of $D_{\K}(R)$ is $\K$.
\end{corollary}
\begin{proof}
	Let $\varphi \in D_{\K}(R)$ be such that $[\varphi , \psi]= 0$ for any 
        $\psi \in D_{\K}(R)$.
In particular, $[\varphi, x_i] = [\varphi , \lambda_{x_i}]=0$ for all 
$i\leq n$.
	That is, $\varphi = \rho_a$ for some $a$ in the centre of $R$. Hence, 
        $a\in \K$.
\end{proof}

\begin{proposition}\label{inner-derivations}
  A derivation of $R$ which is in $D_{\K}^0(R)$ is a sum of inner derivations.
\end{proposition}
\begin{proof}
	  Let $\phi$ be in $D_{\K}^0(R)$ be a derivation.  Since 
          $D_{\K}^0(R)$ is 
          generated
	  by left and right multiplications, we can write 
	    $$\phi = \sum \alpha_i \lambda_{a_i} \rho_{b_i}$$
	  We may assume that if $i \neq j$, then $(a_i, b_i) \neq (a_j,b_j)$.
  
	  Since $\phi$ is a derivation, it has no constant term, so $a_i$ and 
          $b_i$ 
          are not both 1. 

  Let $d = \textup{max}\{\degr{a_i} + \degr{b_i}\}$ and $\tau = x_1^d x_2$.
  Then $\phi(\tau^2) = \sum \alpha_i a_i \tau^2 b_i$.
  Since $\phi$ is a derivation, we also can write
  $\phi(\tau^2) = \sum \alpha_i (a_i \tau  b_i \tau + \tau a_i \tau b_i)$. 
  
  The monomials of this expression have three forms.  
  Put $A_i = a_i \tau^2 b_i$, $B_i = a_i \tau b_i \tau$, and 
  $C_i = \tau a_i \tau b_i$.
  Then $\sum \alpha_i A_i = \sum \alpha (B_i + C_i)$.
  Comparing terms, we have three possible relations among the monomials.

  If $A_i = B_j$, then $\tau b_i = b_j \tau$.  Thus $b_j$ is a left 
  factor of $\tau$
   and $\tau$ is a left factor of $b_j\tau$.
  Thus $b_j = 1$.  It follows that $b_i = 1$ and $a_i = a_j$.  Because the 
  $(a_i, b_i)$
  are chosen to be pairwise distinct, we have $i = j$.
 
  If $A_i = C_j$, then $a_i \tau = \tau a_j$.  Again, we see $a_i = 1$, and so 
  $a_j = 1$, $b_i = b_j$, and $i = j$.
 
  If $B_i = C_j$, then $a_i \tau b_i \tau = \tau a_j \tau b_j$.  
  Thus $a_i = 1$, $b_i = a_j$, and $b_j = 1$.

  In particular, any non-zero term of $\phi(\tau^2)$ must have either 
  $a_i = 1$ or $b_i = 1$.
  Let $I = \{ i \,|\, b_i = 1\}$ and $J = \{ i \,|\, a_i = 1\}$. 
  Since $\phi$ has no constant terms, $I$ and $J$ do not intersect.
  Then 
     $\sum_{i \in I} \alpha_i A_i = \sum_{i \in I} \alpha_i B_i$ and 
     $\sum_{j \in J} \alpha_j A_j = \sum_{j \in J} \alpha_j C_j$.
  Thus $0 = \sum_{i \in I} \alpha_i C_i + \sum_{j \in J} \alpha_j B_j$.
  It follows that for each $i \in I$ there is some $j \in J$ such that
  $\alpha_i = - \alpha_j$, $a_i = b_j$, and $b_i = a_j = 1$.  Let us denote 
  such a $j$ by $j(i)$.

  Thus, $\phi 
  = \sum_I \alpha_i (\lambda_{a_i} - \rho_{b_{j(i)}}) 
  = \sum_I \alpha_i (\lambda_{a_i} - \rho_{a_i})$, 
  and so $\phi$ is a  sum of inner derivations.
\end{proof}

For each $a \in R$, and $i \leq n$ let $\del_i^a$ be the derivation on $R$
defined by $\del_i^a (x_j) = \delta_{i,j} a$. 
\begin{remark}
 Note, for $a \in R$, the inner-derivation 
 $\lambda_a - \rho_a = \sum _{i=1}^n \del_i^{ax_i - x_ia}$.
\end{remark}

\begin{proposition}\label{0-1-order}
	We have the following relations among the operators 
	$\del_i^a , \lambda_b, \rho_b$ for $a,b \in R, i\leq n$.
	\[
	[\del_i^a, \lambda_b]= \lambda_{\del_i^a(b)},
	\quad 	[\del_i^a, \rho_b]= \rho_{\del_i^a(b)}
	\quad [\del_i^a ,\del_j ^b]= \del_j^{\del_i^a (b)}
				- \del_i^{\del_j^b (a)}.
	\]
\end{proposition}
\begin{proof}
	Follows from the definition of $\del_i^a$ for $i\leq n, a \in R$.
\end{proof}

\begin{remark}
One can immediately see that the Lie algebra of derivations on $R$, denoted 
$Der (R)$,  is not simple even  in the 
case when characteristic of $\K$ is $0$.  Let $\bar{R}= \K [t_1,\cdots ,t_n]$ 
denote the polynomial algebra on $n$ variables. The Lie algebra of 
derivations on 
$\bar{R}$, denoted by $Der (\bar{R})$ 
is the vector space spanned by $\{ f\del_{t_i} \mid f\in \bar{R}, 1\leq i 
\leq n \}$ 
where each $\del_{t_i}$ denotes the usual partial derivation with respect to 
$t_i$.  
Consider the quotient algebra map ${}^- : R \to \bar{R}$ given by $\bar{x_i} 
= t_i$ 
for $1\leq i \leq n$. This map gives rise to a map of Lie algebras 
${}^- : Der (R) 
\to Der (\bar{R})$ given by
$\bar{\del_i^a} = \bar{a}\del_{t_i}$. Note that the adjoint derivations are 
in the 
kernel of this map. 

\end{remark}


\begin{definition}\label{Def:usualdef}
For $r=1, I=(i_1)$ and $J=(a_1)$, with $i_1 \leq n, a_1 \in R$, 
set $\del_{I}^{J} = \del_{i_1}^{a_1}$.
For an $r\in \N, r\geq 2$, 
let $I =(i_1, i_2,\cdots ,i_r)$ be a sequence of natural numbers
$i_j \leq n$ and $J=(a_1,\cdots ,a_r)$ be a sequence of elements from $R$.
Further, let $\widehat{I} = (i_2,\cdots ,i_r)$ and
$\widehat{J} = (a_2,\cdots ,a_r)$. 
Denote by $\del_I^J \in D_{\K}^r(R)$ the operator which satisfies the 
commutator rules 
\[
[\del_I^J , x_{i_1}] = a_1\del_{\widehat{I}}^{\widehat{J}}, \quad
[\del_I^J , x_{l}] = 0 \textit{ for }l \neq i_1, \textit{ and }
\quad \del_I^J (1) =0.
\]
\end{definition}
  
\begin{remark}\label{remark-Delta-not-equal-to-D}
Note that $\del_I^J$ is not the same operator as 
$\del_{i_1}^{a_1}\cdots \del_{i_r}^{a_r}$. For instance,
\begin{align*}
[\del_1^{a_1} \del_1^{a_2} , x_1] &= a_1\del_1^{a_2}+ \del_1^{a_1}a_2
	= a_1\del_1^{a_2} + a_2 \del_1^{a_1} + \lambda_{\del_1^{a_1}(a_2)}
\textit{ and}\\
[\del_1^{a_1}\del_1^{a_2}, x_i]&=0 \textit{ for }i\neq 1.
\end{align*}
Therefore,
$\del_1^{a_1} \del_1^{a_2} = \del_{(1,1)}^{(a_1,a_2)} + 
\del_{(1,1)}^{(a_2,a_1)} + 
	\del_1^{\del_1^{a_1}(a_2)}$.
In general,
\[
\del_i^{a_1} \del_j^{a_2} = \del_{(i,j)}^{(a_1,a_2)} +
\del_{(j,i)}^{(a_2,a_1)} + 
	\del_j^{\del_i^{a_1}(a_2)}.
\]
Indeed, for $i=j$, the argument is
 identical to 
the one
given above. For $i\neq j$, note that
\begin{align*}
[\del_i^{a_1}\del_j^{a_2}, x_i] &= a_1\del_j^{a_2}, \\
[\del_i^{a_1}\del_j^{a_2},x_j] &= \del_i^{a_1}a_2 = a_2\del_i^{a_1} + 
\lambda_{\del_i^{a_1}(a_2)},
\textit{ and }\\
[\del_i^{a_1}\del_j^{a_2}, x_k]&=0 \textit{ for }k \neq i,j.
\end{align*}

\noindent
The operator $\del_I^J$ is so defined that $\del_I^J(x_{i_1}\cdots x_{i_r})
=a_1\cdots a_r$ for $I=(x_{i_1},\cdots ,x_{i_r})$ and
$J=(a_1,\cdots ,a_r)$. This can be checked by induction on $r$. 

\noindent
By induction on $r$, we see that $\del_{(i,i,\cdots ,i)}^{(a,a,\cdots ,a)} =
\dfrac{(\del_i^a)^r}{r!}$ for any $a\in R$, and $i\leq n$ when the 
characteristic of $\K$ is $0$.

\noindent
Let $\Delta (R)$ be the associative subalgebra of $D_{\K}(R)$ generated by
$D^1_{\K}(R)$ and 
\[
\Delta ^r(R) = D_{\K}^r(R) \cap \Delta (R) \textit{ for } 
r\geq 1. 
\]
Since $\rho_a = \lambda_a - \sum_{i=1}^n \del_i^{[a,x_i]}$, 
the algebra $\Delta (R)$ is generated by the set 
$\{ \lambda_{x_i}, \del_i^{w} \mid i\leq n, w \textit{ word }\in R \}$.
If $\del^{(1,1)}_{(1,2)} \in \Delta (R)$, then 
$\del^{(1,1)}_{(1,2)} = \lambda_a + \sum f_{i_1,\cdots ,i_r}^{a_1,\cdots ,a_r} 
\del_{i_1}^{a_1}\cdots \del_{i_r}^{a_r}$
for $f_{i_1,\cdots,i_r}^{a_1,\cdots ,a_r}, a \in R$, and words 
$a_1,\cdots , a_r 
\in R$.
Since $\del_{(1,2)}^{(1,1)} (1) =0$, we get $a=0$. Let $c=x_1x_2-x_2x_1$. 
Note that
$\del_{(1,2)}^{(1,1)}(c) = 1$. For any derivation $\varphi$, we see that 
$\varphi (c)
 = [\varphi (x_1),x_2] +
[x_1, \varphi(x_2)]$. In other words, either $\varphi (c)=0$, or is of degree 
greater than 1. 
Thus we get a contradiction.
Therefore $\del^{(1,1)}_{(1,2)} \notin \Delta (R)$. 
Hence $\Delta (R) \neq D_{\K}(R)$. But we have the following theorem. 
\end{remark}
\begin{theorem}\label{symmetric-relation}
For any sequences $I=(i_1,\cdots ,i_r)$
and $J=(a_1,\cdots ,a_r)$ with $1\leq i_s \leq n$ for and $a_s \in R$,  
\[
\sum _{\sigma \in S_r} \del_{\sigma (I)}^{\sigma (J)} 
- \del_{i_1}^{a_1}\del_{i_2}^{a_2}\cdots \del_{i_r}^{a_r} 
\in \Delta ^{r-1}(R);
\]
here $S_r$ denotes the group of permutations over $r$ elements,
$\sigma (I) = (i_{\sigma (1)}, \cdots ,i_{\sigma (r)})$ and
$\sigma (J) = (a_{\sigma (1)}, \cdots ,a_{\sigma (r)})$. 
That is, $\sum _{\sigma \in S_r} \del_{\sigma (I)}^{\sigma (J)} 
\in \Delta ^{r}(R)$.
\end{theorem}
\begin{proof}
From the remark above we see 
\[
\del_{i,j}^{a_1,a_2} + \del_{j,i}^{a_2,a_1} - \del_i^{a_1} \del_j^{a_2} 
	=
	-\del_j^{\del_i^{a_1}(a_2)} \in \Delta^1(R).
\]
For $i\neq i_1,\cdots ,i_r$,
$\left[ \sum _{\sigma \in S_r} \del_{\sigma (I)}^{\sigma (J)}  , x_i \right] 
	=0$
and so  is
$[\del_{i_1}^{a_1}\del_{i_2}^{a_2}\cdots \del_{i_r}^{a_r} ,x_i]=0$.
Now, 
\[
\left[ \sum _{\sigma \in S_r} \del_{\sigma (I)}^{\sigma (J)} , x_{i_t}\right]
= a_{t} \sum _{\sigma \in S_{r-1}} 
	\del_{\sigma (\overline{I})}^{\sigma (\overline{J})}
\]
where $\overline{I}= (i_1, i_2, \cdots , \widehat{i_t},\cdots ,i_r) $ 
and $\overline{J}=(a_1, a_2, \cdots , \widehat{a_{t}},\cdots ,a_r) $, where
$\widehat{b}$ means that $b$ is absent in the sequence. 
Further,
\[
[ \del_{i_1}^{a_1}\del_{i_2}^{a_2}\cdots \del_{i_r}^{a_r}, 
x_{i_t}] = a_t\del_{i_1}^{a_1}\cdots \widehat{\del_{i_t}^{a_t}}\cdots 
	\del_{i_r}^{a_r} + \psi \textit{ for some }\psi \in \Delta^{r-2}(R).
\]
Induction completes the proof.
\end{proof}
\begin{remark}
The algebra $\Delta$ is not simple even when the characteristic of $\K$ is 
$0$. 
Indeed, let $\mathcal{I}$ be the commutator ideal of $R$;
that is,  let $\mathcal{I}$ be the two sided ideal of $R$ generated by the set
$\{ [a,b] \mid a,b\in R \}$.  If we let $\overline{a}$ denote the image of 
$a$ in
$R/\mathcal{I}$ then $R/\mathcal{I}$ is the polynomial algebra in $n$ 
variables,
$\overline{x_1},\cdots ,\overline{x_n}$. Note, for $a\in R$ and $i\leq n$, 
we have
\[
\lambda_a (\mathcal{I})\subset \mathcal{I}, \quad \rho_a (\mathcal{I}) \subset 
\mathcal{I},
\quad \textit{ and } \quad \del_i^a(\mathcal{I}) \subset \mathcal{I}.
\]
Hence, $\Delta (\mathcal{I}) \subset \mathcal{I}$. We have an algebra 
homomorphism
\begin{align*}
\Delta &\to D(R/\mathcal{I}) \textit{ given by }\\
\lambda_a &\mapsto \lambda_{\overline{a}}, \\
\rho_a &\mapsto \rho_{\overline{a}} =\lambda_{\overline{a}}, \\
\del_i^a &\mapsto \overline{a}\del_i \quad (\textit{here}, \del_i 
\textit{ denotes 
  the partial derivative on polynomial
	algebra}).
\end{align*}
Since $D_{\K}(R/\mathcal{I})$ is generated as an algebra by 
$D_{\K}^1(R/\mathcal{I})$, this map is
surjective with nontrivial kernel (for example, $\lambda_a - \rho_a$ is in the 
kernel for every
$a\in R$).

Moreover, the set $\{ \del_{i_1}^{a_1}\del_{i_2}^{a_2}\cdots \del_{i_r}^{a_r} 
\}_
{r\geq 0}$ does not
form a free basis of $\Delta (R)$ over $R$ as in the case of polynomial 
algebra, 
because 
$\del_i^a \del_j^b- \del_j^b\del_i^a - \del_j^{\del_i^b(a)} + 
\del_i^{\del_i^a(b)} 
=0$. 
Even if we consider the subset $\{ \del_{i_1}^{a_1}\del_{i_2}^{a_2}\cdots 
\del_{i_r}^{a_r} \}_{r\geq 0, i_1\leq i_2 \leq 
\cdots \leq i_r}$  we still do not have a free basis because
$\del_i^a \del_i^b- \del_i^b\del_i^a - \del_i^{\del_i^b(a)} + 
\del_i^{\del_i^a(b)} 
=0$. 

The following two formulae can be checked by induction on $r$.
\begin{align*}
\del_{i_1}^{a_1}\del_{i_2}^{a_2}\cdots \del_{i_r}^{a_r} \lambda_b &= 
	\sum_{s=0}^r \sum_{j_1<j_2<\cdots <j_s} 
        \lambda_{ \del_{i_{j_1}}^{a_{j_1}} 
          \del_{i_{j_2}}^{a_{j_2}}
	\cdots  \del_{i_{j_s}}^{a_{j_s}} (b) }\prod _{t=1, t\notin 
        \{j_1,\cdots ,j_s\}}^r \del_{i_t}^{a_t};\\
\del_{i_1}^{a_1}\del_{i_2}^{a_2}\cdots \del_{i_r}^{a_r} \rho_b &= 
	\sum_{s=0}^r \sum_{j_1<j_2<\cdots <j_s} 
        \rho_{ \del_{i_{j_1}}^{a_{j_1}} 
          \del_{i_{j_2}}^{a_{j_2}}
	\cdots  \del_{i_{j_s}}^{a_{j_s}} (b) }\prod _{t=1, t\notin 
        \{j_1,\cdots ,j_s\}}^r \del_{i_t}^{a_t}.
\end{align*}
\end{remark}
\begin{proposition}\label{commutator-properties}
Here are some properties of the operators $\del_I^J$:
\begin{enumerate}
\item	For $r\geq 1$
	\[
	\del_{(i_1,\cdots ,i_r)}^{(a_1,\cdots ,a_r)}(x_{t_1}\cdots x_{t_r})
	= 	\begin{cases}
		a_1\cdots a_r &\textit{ if } 
			(t_1, \cdots ,t_r) = (i_1,\cdots ,i_r)\\
		0 &\textit{ if } (t_1, \cdots ,t_r) \neq (i_1,\cdots ,i_r)
         	\end{cases}
	\]
	Further, 
	$\del_{(i_1,\cdots ,i_r)}^{(a_1,\cdots ,a_r)}(x_{t_1}\cdots x_{t_k})
	=0$ for $k<r$.
\item
	$
	[\del_{(i_1,\cdots ,i_r)}^{(a_1,\cdots ,a_r)} , \lambda_a] 
	= \sum_{j=2}^r 
		\lambda_{\del_{(i_1,\cdots ,i_{j-1})}^{(a_1,\cdots ,a_{j-1})} 
		(a)}
			\del_{(i_j,\cdots ,i_r)}^{(a_j,\cdots ,a_r)}
	+ \lambda_{\del_{(i_1,\cdots ,i_r)}^{(a_1,\cdots ,a_r)} (a)}
	$ for $r\geq 2$. 

	For $r=1$, 
	$[\del_i^a , \lambda_b] = \lambda_{\del_i^a(b)}$.
	In particular, 
	$[\del_{(i_1,\cdots ,i_r)}^{(a_1,\cdots ,a_r)} , \lambda_a] 
	\in D^{r-1}_{\K}(R)$.
\item
	$
	[\del_{(i_1,\cdots ,i_r)}^{(a_1,\cdots ,a_r)} , \rho_a] 
	= \sum_{j=2}^r \rho_{\del_{(i_j,\cdots ,i_r)}^{(a_j,\cdots ,a_r)} (a)}
			\del_{(i_1,\cdots ,i_{j-1})}^{(a_1,\cdots ,a_{j-1})}
	+ \rho_{\del_{(i_1,\cdots ,i_r)}^{(a_1,\cdots ,a_r)} (a)}
	$ for $r\geq 2$. 

	For $r=1$, 
	$[\del_i^a , \rho_b] = \rho_{\del_i^a(b)}$.
	In particular,
		$
	[\del_{(i_1,\cdots ,i_r)}^{(a_1,\cdots ,a_r)} , \rho_a] 
	\in D^{r-1}_{\K}(R)$.
\end{enumerate}
\end{proposition}
\begin{proof}
\begin{enumerate}
\item	Can be checked by induction on $r$.
\item	For $r=1$, the conclusion 
	is in part (2) of proposition \ref{0-1-order}. For $r>1$ we proceed
	by induction, assuming the result for earlier cases. Let
	$I= (i_1,\cdots ,i_r)$, $J=(a_1,\cdots ,a_r)$ and
	$a=x_lb$. Assume the claim for $[\del_I^J, \lambda _b]$. We now
	prove the claim for $[\del_I^J, \lambda _a]$.
	Consider two cases:
	
	\noindent
	Case $i_1\neq l$. Here, 
	\begin{align*}
	\del_{I}^{J} \lambda_{x_l} \lambda_b
		&= \lambda_{x_l}( 
			\del_{I}^{J}\lambda_b)
						\\
		&= \lambda_{x_l} \left( \lambda_b \del_I^J
		+
		\sum_{j=2}^r 
		\lambda_{\del_{(i_1,\cdots ,i_{j-1})}^{(a_1,\cdots ,a_{j-1)}} 
				(b)}
			\del_{(i_j,\cdots ,i_r)}^{(a_j,\cdots ,a_r)}
	+ \lambda_{\del_{I}^{J} (b)} \right)
	\end{align*}
	which gives the claim.	

	\noindent
	Case $i_1 = l$. Here,
	\begin{align*}
	\del_I^J\lambda_{x_l} \lambda_b 
		=& \lambda_{x_l} \del_I^J\lambda_b
			+ \lambda_{a_1}
			\del_{(i_2,\cdots ,i_r)}^{(a_2,\cdots ,a_r)} 
				\lambda_b \\
		=& \lambda_{x_l}\left( 
		\lambda_b \del_I^J + \sum_{j=2}^r 
			\lambda_{\del_{(i_1,\cdots ,i_{j-1})}
			^{(a_1,\cdots ,a_{j-1})} (b)} 
			\del_{(i_j,\cdots ,i_r)}^{(a_j,\cdots, a_r)}
			+ \lambda_{\del_I^J(b)} \right) +
          \end{align*}
          \begin{align*}
		&+ \lambda_{a_1}\left( \lambda_b \del_{(i_2,\cdots ,i_r)}^
				{(a_2,\cdots ,a_r)} +
		\sum_{j=3}^r 
			\lambda_{\del_{(i_2,\cdots ,i_{j-1})}
				^{(a_2,\cdots ,a_{j-1})} 
		(b)}
			\del_{(i_j,\cdots ,i_r)}^{(a_j,\cdots ,a_r)}
		+ \lambda_{\del_{(i_2,\cdots ,i_r)}^{(a_2,\cdots ,a_r)} (b)}
		\right) \\
		=& \lambda_{x_lb}\del_I^J + 
		\sum_{j=3}^r \lambda_{(x_l\del_{(i_1,\cdots ,i_{j-1})}
			^{(a_1,\cdots ,a_{j-1})} +a_1
			\del_{(i_2,\cdots ,i_{j-1})}
			^{(a_2,\cdots ,a_{j-1})} )(b)} 
			\del_{(i_j,\cdots ,i_r)}
			^{(a_j,\cdots ,a_r)} \\
		&+ \lambda_{(x_l\del_{i_1}^{a_1}+a_1) (b)}
			\del_{(i_2,\cdots ,i_r)}
			^{(a_2,\cdots ,a_r)}
		 + \lambda_{x_l\del_{(i_1,\cdots ,i_r)}
			^{(a_1,\cdots ,a_r)}(b) + a_1
			\del_{(i_2,\cdots ,i_r)}
			^{(a_2,\cdots ,a_r)}(b) }\\
		=& \lambda_a \del_I^J +
		\sum_{j=2}^r \lambda_{\del_{(i_1,\cdots ,i_{j-1})}
			^{(a_1,\cdots ,a_{j-1})} (a)} 
			\del_{(i_j,\cdots ,i_r)}
			^{(a_j,\cdots ,a_r)}
			+\lambda_{\del_I^J(a)}.
	\end{align*}
\item	We address this by induction on $r$. For $r=1$, the conclusion 
	is in part (2) of proposition \ref{0-1-order}.
	For $r\geq 2$, consider
	\begin{align*}
	[\del_{(i_1,\cdots ,i_r)}^{(a_1,\cdots ,a_r)}\rho_a 
	- \rho_{\del_{(i_1,\cdots ,i_r)}^{(a_1,\cdots ,a_r)}(a)} 
	, x_{l}]
	&= \delta_{i_1, l} \del_{(i_2,\cdots ,i_r)}^{(a_2,\cdots ,a_r)}\rho_a 
	\textit{ and }\\
	\left( 
	\del_{(i_1,\cdots ,i_r)}^{(a_1,\cdots ,a_r)}\rho_a 
	- \rho_{\del_{(i_1,\cdots ,i_r)}^{(a_1,\cdots ,a_r)}(a)} 
	\right) (1) &= 0.
	\end{align*}
	Now, by induction assumption,
	\[
	\del_{(i_2,\cdots ,i_r)}^{(a_2,\cdots ,a_r)}\rho_a 
	= 
	\rho_a 	\del_{(i_2,\cdots ,i_r)}^{(a_2,\cdots ,a_r)} +
	\sum_{j=3}^r \rho_{\del_{(i_j,\cdots ,i_r)}^{(a_j,\cdots ,a_r)} (a)}
			\del_{(i_2,\cdots ,i_{j-1})}^{(a_2,\cdots ,a_{j-1})}
	+ \rho_{\del_{(i_2,\cdots ,i_r)}^{(a_2,\cdots ,a_r)} (a)} .
	\]
	Therefore,
	\begin{align*}
	\del_{(i_1,\cdots ,i_r)}^{(a_1,\cdots ,a_r)}\rho_a 
	- \rho_{\del_{(i_1,\cdots ,i_r)}^{(a_1,\cdots ,a_r)}(a)} 
	=& 
	\rho_a 	\del_{(i_1,\cdots ,i_r)}^{(a_1,\cdots ,a_r)}
        +
	\sum_{j=3}^r \rho_{\del_{(i_j,\cdots ,i_r)}^{(a_j,\cdots ,a_r)} (a)}
			\del_{(i_1,\cdots ,i_{j-1})}^{(a_1,\cdots ,a_{j-1})}\\
	&+ \rho_{\del_{(i_2,\cdots ,i_r)}^{(a_2,\cdots ,a_r)} (a)} 
	\del_{i_1}^{a_1}
	\end{align*}
	which gives us the claim.
\end{enumerate}
\end{proof}
\begin{remark}
For $I=(i_1,\cdots ,i_r)$, $A=(a_1,\cdots ,a_r)$  and 
$J = (i_l,\cdots ,i_{l+t})$ 
for some $1\leq l \leq l+t\leq r$, let 
$A_J =(a_l,\cdots ,a_{l+t})$. For two finite sequences $J,K$, by $(J,K)$ 
we mean
the concatenation of $J$ and $K$.
Then Part (1) or (2) of the above proposition gives
\[
\del_I^A (ab) = \sum _{(J,K)=I} \del_J^{A_J}(a) \del_K^{A_K} (b).
\]
Part (1) of the above Proposition has an immediate important corollary.
\end{remark}
\begin{corollary}
$R$ as a left $D_{\K}(R)$ module is simple.
\end{corollary}
\begin{proof}
Let any $r\in R, r\neq 0, r\notin \K$. We need to show the existence of 
$\Phi \in D_{\K}(R)$ such that $\Phi (r) \in \K^*$.
Suppose degree of $r$ is $n$, $n\geq 1$. Let 
$\alpha_{i_1,\cdots ,i_l}x_{i_1}^{k_1}x_{i_2}^{k_2}\cdots x_{i_l}^{k_l}$ be a  
summand of $r$ of degree $n$. That is,
\[
r = \alpha_{i_1,\cdots ,i_l}x_{i_1}^{k_1}x_{i_2}^{k_2}\cdots x_{i_l}^{k_l}
	+ \textit{terms of degree less than or equal to }n.
\]

Let $I=(\underbrace{i_1,\cdots ,i_1}_{k_1 \textit{ times }},
\underbrace{i_2,\cdots ,i_2}_{k_2 \textit{ times }},
\cdots , \underbrace{i_l,\cdots ,i_l}_{k_l \textit{ times }})$ and
$J = ( \underbrace{1,1,\cdots ,1}_{n \textit{ times }})$.
Then $\del_I^J (r) = \alpha_{i_1,\cdots ,i_l} \in \K^*$.
\end{proof}

\newcommand{\PermGroup}[1]{S_{#1}}
\newcommand{\ShufflingSet}[2]{T_{({#1},{#2})}}
\newcommand{\permutationExtension}[2]{{#1}^{({#2})}}
\newcommand{\sigmak}[1]{\permutationExtension{\sigma}{#1}}
\newcommand{\tauk}[1]{\permutationExtension{\tau}{#1}}

Let $\PermGroup{r+s}$ be the group of permutations of $r+s$ elements, and let 
$\ShufflingSet{r}{s}$ be the subset
of $\PermGroup{r+s}$ which preserves the order of the first $r$ elements and 
the 
last $s$ elements.  That is, $\tau \in \ShufflingSet{r}{s}$ is increasing on 
$1, 
\ldots, r$ and on $r+1,\ldots, r+s$.  Note that $\ShufflingSet{0}{s}$ and 
$\ShufflingSet{r}{0}$ consist of just the identity permutation.

For any $\sigma \in \PermGroup{r+s}$ and $k$, $1 \leq k \leq r+s$, we can 
define a 
permutation $\sigmak{k} \in \PermGroup{r+s+1}$ such that
$$
\sigmak{k}(i) = 
  \left\{
    \begin{array}{ll}
          \sigma (i)+1 & i < k\\
	1 & i = k \\
      	\sigma(i-1)+1 & i > k
    \end{array}
  \right.
 $$
That is, $\sigmak{k}$ takes $k$ to $1$ and acts like $\sigma$ on the remaining 
symbols.
As $\sigma$ is a bijection from $\{1, \ldots, r+s\}$ to $\{1, \ldots, r+s \}$,
 we 
can see that $\sigmak{k}$ is a bijection from $\{1, \ldots, \widehat{k}, 
\ldots, 
r+s+1\}$ to $\{2, \ldots, r+s+1\}$. Hence $\sigmak{k}$ is indeed a permutation.

\begin{lemma}
  $\ShufflingSet{r}{s} =
    \lbrace \tauk{1} \,|\, \tau \in \ShufflingSet{r-1}{s} \rbrace \cup
    \lbrace \tauk{r+1} \,|\, \tau \in \ShufflingSet{r}{s-1} \rbrace$.
\end{lemma}
\begin{proof}
  First, let $\tau \in \ShufflingSet{r-1}{s}$.  We will show $\tauk{1} \in 
  \ShufflingSet{r}{s}$.

  As $\tauk{1}$ on $2, \ldots, r+s$ is just a translation of $\tau$, and 
  $\tau$ is
  increasing on $1, \ldots, r-1$ and on $r, \ldots, r-1 + s$, we have 
  $\tauk{1}$ is
  increasing on $2, \ldots, r$ and on $r+1, \ldots, r+s$ respectively.  Since 
  $\tauk{1}(1) = 1$, 
  we see $\tauk{1}\in \ShufflingSet{r}{s}$.

  Second, let $\tau \in T(r, s-1)$.  We will show $\tauk{r+1} \in 
  \ShufflingSet{r}{s}$.
  As $\tau$ is increasing on $1,\ldots,r$, so is $\tauk{r+1}$.
  As $\tau$ is increasing on $r+1, \ldots, r+s-1$, so $\tauk{r+1}$ is
  increasing on $r+2, \ldots, r+s$.  As $\tauk{r+1}(r+1) = 1 
  < \tauk{r+1}(r+2)$, 
  we have $\tauk{r+1}\in \ShufflingSet{r}{s}$.
  
  This proves $\ShufflingSet{r}{s} \supseteq \lbrace \tauk{1} \,|\, \tau \in 
  \ShufflingSet{r-1}{s} \rbrace \cup
    \lbrace \tauk{r+1} \,|\, \tau \in \ShufflingSet{r}{s-1} \rbrace$.
  
  Now, let $\tau \in \ShufflingSet{r}{s}$.  Then either $\tau(1) = 1$ or 
  $\tau(r+1) = 1$.
 
  Case $\tau(1) = 1$.  Let $\sigma \in \PermGroup{r+s-1}$ be such that
   $\sigma(i)=\tau(i+1)-1$.
  Then $\sigmak{1} = \tau$.
  As $\tau$ is increasing on $2,\ldots, r$ and on $r+1,\ldots, r+s$, 
  we have $\sigma$ is increasing on $1,\ldots,r-1$ and on 
  $r, \ldots, r -1 +s$, 
  respectively. 
  Thus $\sigma \in \ShufflingSet{r-1}{s}$.

  Case $\tau(r+1) = 1$. Let $\sigma \in \PermGroup{r+s-1}$ be such that
   $\sigma(i)=\tau(i)-1$ for $1 \leq i \leq r$ and
   $\sigma(i)=\tau(i+1)-1$ for $r+1\leq i \leq r+s-1$.
  Then $\sigmak{r+1} = \tau$.

  As $\tau$ is increasing on $1, \ldots, r$ and on $r+2, \ldots, r+s$, 
  we have $\sigma$ is increasing on $1, \ldots, r$ and on $r+1, \ldots, r+s -1$
  respectively.  Thus $\sigma \in \ShufflingSet{r}{s-1}$.
\end{proof} 
Recall that
$(I,J)$ denotes the concatenation  of $I$ with $J$, and
$(A,B)$ denotes the concatenation of $A$ with $B$.
\begin{theorem}\label{theorem-product of dels}
  $$\del_I^A \del_J^B =
    \sum_{\tau \in \ShufflingSet{r}{s}} \del_{\tau((I,J))}^{\tau((A,B))} +
    \sum \lbrace \del_K^C \,|\, \text{$\del_K^C$ is of order less than $r+s$}.
    \rbrace$$
\end{theorem}

\begin{proof}
  Put $N=r+s$. We proceed by induction on $N$.
  
  Our base case starts with $N=2$ and $r=s=1$.  
  Corollary \ref{remark-Delta-not-equal-to-D}
   proves the result. 

  \newcommand{\Ihat}{\widehat{I}}
  \newcommand{\Jhat}{\widehat{J}}
  \newcommand{\Ahat}{\widehat{A}}
  \newcommand{\Bhat}{\widehat{B}}
  \newcommand{\dA}{\del_{I}^{A}}
  \newcommand{\dB}{\del_{J}^{B}}

  \newcommand{\dAhat}{\del_{\Ihat}^{\Ahat}}
  \newcommand{\dBhat}{\del_{\Jhat}^{\Bhat}}
  \newcommand{\dAhatB}{\del_{\IhatJ}^{\Ahat B}}
  \newcommand{\dABhat}{\del_{I\Jhat}^{A \Bhat}}

  Suppose that the theorem holds for $N < N_0$ and that $r+s=N_0$.  Then 
  $$\dAhat \dB = \sum_{\tau \in \ShufflingSet{r-1}{s}} 
  \del_{\tau((\Ihat , J))}^{\tau((\Ahat ,B))} + \sum_{\alpha} 
  \del_{L_{\alpha}}^{C_{\alpha}}$$
  $$\dA \dBhat = \sum_{\tau \in \ShufflingSet{r}{s-1}} 
  \del_{\tau((I, \Jhat))}^{\tau((A, \Bhat))} + \sum_{\beta} 
  \del_{L_{\beta}}^{C_{\beta}}$$
  $$\dAhat \dBhat = \sum_{\tau \in \ShufflingSet{r-1}{s-1}} 
  \del_{\tau((\Ihat, \Jhat))}^{\tau((\Ahat, \Bhat))} + \sum_{\gamma} 
  \del_{L_{\gamma}}^{C_{\gamma}}$$
  where $\del_{L_{\alpha}}^{C_{\alpha}}$, $\del_{L_{\beta}}^{C_{\beta}}$, and 
  $\del_{L_{\gamma}}^{C_{\gamma}}$
  are all of order $N_0-2$ or lower.

\newcommand{\bracket}[2]{\left[{#1},{#2}\right]}
\newcommand{\dell}[2]{\del_{#1}^{#2}}
  For each $i$, put $f_i=\dell{i_1}{a_1}(x_i) = \delta_{ii_1}a_1$, 
  $g_i=\dell{j_1}{b_1}(x_i)=\delta_{ij_1}b_1$, and $h_i=\dell{i_1}{a_1}(g_i)=
  \delta_{ij_1}\dell{i_1}{a_1}(b_1)$.  Then we have 
  \begin{align*}
    \bracket{\dA\dB}{\lambda_{x_i}}
     =& \bracket{\dA}{\lambda_{x_i}} \dB + \dA \bracket{\dB}{\lambda_{x_i}}\\
     =& \lambda_{f_i} \dAhat \dB + \dA \lambda_{g_i} \dBhat \\
     =& \lambda_{f_i} \dAhat \dB + \lambda_{g_i}\dA \dBhat + \lambda_{h_i} 
     \dAhat 
     \dBhat \\
     =& \sum_{\tau \in \ShufflingSet{r-1}{s}}\lambda_{f_i} 
     \dell{\tau((\Ihat, J))}{\tau((\Ahat, B))}
       + \sum_{\tau \in \ShufflingSet{r}{s-1}} \lambda_{g_i} 
       \dell{\tau((I, \Jhat))}{\tau((A, \Bhat))} \\
      & + \sum_{\tau \in \ShufflingSet{r-1}{s-1}} \lambda_{h_i}
       \dell{\tau((\Ihat, \Jhat))}{\tau((\Ahat, \Bhat))} 
       + \sum_{\alpha} \lambda_{f_i} \dell{L_{\alpha}}{C_{\alpha}} 
       + \sum_{\beta} \lambda_{g_i}\dell{L_{\beta}}{C_{\beta}}
       + \sum_{\gamma} \lambda_{h_i} \dell{L_{\gamma}}{C_{\gamma}}\\
  =& \sum_{\tau \in \ShufflingSet{r-1}{s}} 
  \bracket{\dell{(i_1,\tau((\Ihat, J)))}{(a_1,\tau((\Ahat, B)))}}
  {\lambda_{x_i}}
      + \sum_{\tau \in \ShufflingSet{r}{s-1}} 
      \bracket{\dell{(j_1,\tau((I, \Jhat)))}{(b_1,\tau((A, \Bhat)))}}
      {\lambda_{x_i}}\\ 
     & + \sum_{\tau \in \ShufflingSet{r-1}{s-1}} 
      \bracket{\dell{(j_1,\tau(\Ihat, \Jhat))}
        {(\dell{i_1}{a_1}(b_1),\tau((\Ahat, \Bhat)))}}{\lambda_{x_i}} 
   + \sum_{\alpha}   \bracket{\dell{(i_1,L_{\alpha})}
     {(a_1,C_{\alpha})}}{\lambda_{x_i}}\\
    & + \sum_{\beta}     \bracket{\dell{(j_1, L_{\beta})}
       {(b_1, C_{\beta})}}{\lambda_{x_i}}
     + \sum_{\gamma} \bracket{\dell{(j_1, L_{\gamma})}
       {(\dell{i_1}{a_1}(b_1), C_{\gamma})}}{\lambda_{x_i}}
  \end{align*}
Note that
$(i_1,\tau((\Ihat, J))) = \tauk{1}((I, J))$, $(j_1,\tau((I \Jhat)) )= 
\tauk{r+1}((I ,J))$, \\
$(a_1, \tau((\Ahat, B))) = \tauk{1}((A,B))$, and $(b_1, \tau((A, \Bhat))) 
= \tauk{r+1}((A,B))$.
Thus, by our Lemma, 
$$\sum_{\tau \in \ShufflingSet{r}{s}} \dell{\tau((I,J))}{\tau((A,B))} 
=\sum_{\tau \in \ShufflingSet{r-1}{s}}  \dell{(i_1,\tau((\Ihat, J)))}
{(a_1,\tau((\Ahat, B)))}
+\sum_{\tau \in \ShufflingSet{r}{s-1}} \dell{(j_1,\tau((I, \Jhat)))}
{(b_1,\tau((A, \Bhat)))}.$$

Put 
\begin{align*}
\phi
 =& \sum_{\tau \in \ShufflingSet{r}{s}} \dell{\tau((I, J))}{\tau((A, B))}
+  \sum_{\tau \in \ShufflingSet{r-1}{s-1}} \dell{(j_1\tau((\Ihat, \Jhat)))}
{(\dell{i_1}{a_1}(b_1),\tau((\Ahat, \Bhat)))}
+ \sum_{\alpha}    \dell{(i_1,L_{\alpha})}{(a_1,C_{\alpha})}\\
&+ \sum_{\beta}    \dell{(j_1, L_{\beta})}{(b_1, C_{\beta})}
+ \sum_{\gamma} \dell{(j_1, L_{\gamma})}{(\dell{i_1}{a_1}(b_1), C_{\gamma})}.
\end{align*}
Since $\phi(1) = \dA \dB(1) = 0$ and $\bracket{\phi}{\lambda_{x_i}} = 
\bracket{\dA\dB}{\lambda_{x_i}}$ for each $i$, we have $\phi = \dA\dB$.
Regardless of how the tuples are permuted, $\dell{(j_1,\Ihat ,\Jhat)}
{(\dell{i_1}{a_1}(b_1),\Ahat, \Bhat)}$ has order $N_0-1$.  Likewise, each
$\dell{(i_1,L_{\alpha})}{(a_1,C_{\alpha})}$, $\dell{(j_1, L_{\beta})}
{(b_1, C_{\beta})}$, and 
$\dell{(j_1, L_{\gamma})}{(\dell{i_1}{a_1}(b_1), C_{\gamma})}$ 
are of order $N_0-1$ or less.  Thus our expression $\phi$ is of the required 
form.
\end{proof}
\begin{remark}
The lower order terms which appear in the statement of Theorem 
\ref{theorem-product of dels} can be described. 
If we denote $\del_{(i_1,\cdots ,i_r)}^{(a_1,\cdots ,a_r)} $ by 
$\begin{pmatrix}a_1, &a_2, &\cdots , &a_r\\  i_1, &i_2, &\cdots ,&i_r  
                     \end{pmatrix}$
then every lower order term looks like
\[
\begin{pmatrix}
a_1,&\cdots ,&a_{p_1-1}, &d_1(b_1),
&a_{p_{q_1}+1}, &\cdots ,&a_{p_{q_2}-1}, 
 &d_2(b_2), &\cdots \\
i_1,&\cdots ,&i_{p_1-1}, &j_1,
&a_{p_{q_1}+1}, &\cdots ,&a_{p_{q_2}-1}, 
 &j_2, &\cdots 
\end{pmatrix}
\]
where $d_1(b_1) = \del_{(i_{p_1},i_{p_1+1},\cdots , 
  i_{p_{q_1}})}^{(a_{p_1},a_{p_1+1},\cdots , a_{p_{q_1}})}(b_1)$
and $d_2(b_2)= \del_{(i_{p_2},i_{p_2+1},\cdots , i_{p_{q_2}})}^
 {(a_{p_2},a_{p_2+1},\cdots , a_{p_{q_2}})}(b_2)$.

The last term is
\[
\begin{pmatrix}
b_1,&b_2,&\cdots, &b_{s-1}, &\del_I^A (b_s)\\
j_1, &j_2, &\cdots, &j_{s-1}, &j_s
\end{pmatrix}
.
\]
In particular, the operators which appear in the product $\del_I^A\del_J^B$ 
are of order at least $s$, which is the order of
$\del_J^B$. For proof, adapt the proof of Theorem 
\ref{theorem-product of dels} 
along with part (2) of 
Proposition \ref{commutator-properties}.
We present a few examples.
\begin{align*}
\begin{pmatrix}
a\\
i
\end{pmatrix}
\cdot
\begin{pmatrix}
b_1, &\cdots, &b_s\\
j_1, &\cdots, &j_s
\end{pmatrix}
=&
\sum_{l=1}^s 
\begin{pmatrix}
b_1,&\cdots ,&b_{l-1}, &a,&b_l, &\cdots, &b_s\\
j_1, &\cdots, &j_{l-1},&i,&j_l,&\cdots, &j_s
\end{pmatrix}\\
+&
\sum_{l=1}^s
\begin{pmatrix}
b_1, &\cdots ,&b_{l-1}, &\del_i^a(b_l),&b_{l+1}, &\cdots ,&b_s\\
j_1,&\cdots, &j_{l-1}, &j_l,&j_{l+1},&\cdots, &j_s
\end{pmatrix}.
\end{align*}
\begin{align*}
&
\begin{pmatrix}
a_1,&\cdots, &a_r\\
i_1,&\cdots, &i_r
\end{pmatrix}
\cdot
\begin{pmatrix}
b\\
j
\end{pmatrix}
=
\sum_{l=1}^r 
\begin{pmatrix}
a_1,&\cdots, &a_{l-1}, &b, &a_l,&\cdots, &a_r\\
i_1,&\cdots, &i_{l-1}, &j,&j_l, &\cdots, &i_r
\end{pmatrix}\\
+&
\sum_{p=1}^r \sum_{s=0}^{r-p}
\begin{pmatrix}
a_1,&\cdots, &a_{p-1}, &\del_{(i_p,\cdots ,i_{p+s})}^{(a_p,\cdots, a_{p+s})} 
(b), 
&a_{p+s+1}, &\cdots, &a_r\\
i_1,&\cdots ,&i_{p-1},&j, &i_{p+s+1}, &\cdots, &i_r
\end{pmatrix}.
\end{align*}
\begin{align*}
\begin{pmatrix}
a_1, &a_2\\
i_1, &i_2 
\end{pmatrix}
\cdot
\begin{pmatrix}
b_1, &b_2\\
j_1, &j_2 
\end{pmatrix}
=& 
\begin{pmatrix}
a_1, &a_2, &b_1, &b_2\\
i_1, &i_2, &j_1, &j_2 
\end{pmatrix}
+
\begin{pmatrix}
a_1, &b_1, &a_2, &b_2\\
i_1, &j_1, &i_2, &j_2 
\end{pmatrix}
\\
+&
\begin{pmatrix}
b_1, &a_1, &a_2, &b_2\\
j_1, &i_1, &i_2, &j_2 
\end{pmatrix}
+
\begin{pmatrix}
a_1, &b_1, &b_2, &a_2\\
i_1, &j_1, &j_2, &i_2 
\end{pmatrix}
\\
+&
\begin{pmatrix}
b_1, &a_1, &b_2, &a_2\\
j_1, &i_1, &j_2, &i_2 
\end{pmatrix}
+
\begin{pmatrix}
b_1, &b_2, &a_1, &a_2\\
j_1, &j_2, &i_1, &i_2 
\end{pmatrix}
\\
+&
\begin{pmatrix}
a_1,  &\del_{i_2}^{a_2}(b_1), &b_2\\
i_1,  &j_1, &j_2 
\end{pmatrix}
+
\begin{pmatrix}
a_1, &b_1, &\del_{i_2}^{a_2}(b_2)\\
i_1,  &j_1, &j_2 
\end{pmatrix}
\\
+&
\begin{pmatrix}
b_1, &a_1, &\del_{i_2}^{a_2}(b_2)\\
j_1,  &i_1, &j_2 
\end{pmatrix}
+
\begin{pmatrix}
\del_{i_1}^{a_1}(b_1),  &{a_2}, &b_2\\
j_1,  &i_2, &j_2 
\end{pmatrix}
\\
+&
\begin{pmatrix}
\del_{i_1}^{a_1}(b_1),  &{b_2}, &a_2\\
j_1,  &j_2, &i_2 
\end{pmatrix}
+
\begin{pmatrix}
b_1, &\del_{i_1}^{a_1}(b_2),  &a_2\\
j_1,  &j_2, &i_2 
\end{pmatrix}
\\
+&
\begin{pmatrix}
\del_{i_1}^{a_1}(b_1),  &\del_{i_2}^{a_2}(b_2)\\
j_1,  &j_2
\end{pmatrix}
+
\begin{pmatrix}
b_1,  &\del_{i_1,i_2}^{a_1,a_2}(b_2)\\
j_1,  &j_2
\end{pmatrix}.
\end{align*}
\end{remark}
\begin{corollary}\label{commuting dels}
	$
	[\del_I^J, \del_K^L] \in D^{r+t-1}_{\K}(R) 
	$ for 
 	$\del_I^J \in D^r_{\K}(R)$, $\del_K^L \in D^t_{\K}(R)$,
	and $r,t\geq 1$.
\end{corollary}

\begin{theorem}\label{gen-D(free)}
For $r\geq 1$ $D^r_{\K}(R)$ is generated as a left $D^0_{\K}(R)$-module by 
the set
\[
\{ \del_I^J \mid I =(i_1,\cdots ,i_s), J =(a_1,\cdots ,a_s), 
 i_j\in \N, i_j\leq n, a_j \in R, 1\leq s\leq r \} 
\cup 
\{ \lambda_1\}
\]
\end{theorem}
\begin{proof}
	Let $\varphi \in D^1_{\K}( R)$ such that 
	$[ \varphi ,x_i] =\sum_j \lambda_{a_{i,j}} 
	\rho_{b_{i,j}} \in D^0_{\K}(R)$,
	for $a_{i,j},b_{i,j}\in R$ and $i\leq n$. Then
	$\psi = \varphi - \sum_i \sum_j \rho_{b_{i,j}}\del_i^{a_{i,j}}$ 
	is such that $[\psi ,x_i]=0$ for every $i\leq n$. Therefore 
	$\psi = \rho_r$ for some $r\in R$. Hence,
	$\varphi = \rho_r +\sum_{i,j} \rho_{b_{i,j}}\del_i^{a_{i,j}}$.

	In general, let
	$\varphi \in D^r_{\K}( R)$ such that 
	$[ \varphi ,x_i] =\sum_{i,I,J} \rho _{a^i_{I,J}} \lambda_{b^i_{I,J}}
	\del_I^J \in D^{r-1}_{\K}(R)$,
	for $a^i_{I,J},b^i_{I,J} \in R$, $i\leq n$ 
	and $I,J$ the appropriate
	sequences of length less than or equal to $r-1$. Then
	$\psi = \varphi - \sum_{i,I,J} \rho_{a^i_{I,J}} 
			\del_{(i,I)}^{(b^i_{I,J}, J)}$
	is such that $[\psi ,x_i]=0$ for every $i\leq n$. Therefore 
	$\psi = \rho_s$ for some $s\in R$. Hence,
	$\varphi = \rho_s +\sum_{i,I,J} \rho_{a^i_{I,J}} 
			\del_{(i,I)}^{(b^i_{I,J}, J)}$;
	here, $(i,I)$ is the concatenation of
	$(i)$ and $I$, while $(b^i_{I,J}, J)$ is the concatenation
	of $ (b^i_{I,J})$ and $J$.
\end{proof}

Suppose $J=(a_1,\cdots ,a_s)$ is such that
$a_p = \alpha_1 b_{p_1} + \alpha_2 b_{p_2}$ for some $p\leq s$, 
$\alpha_1, \alpha_2 \in \K$, and $b_{p_1}, b_{p_2} \in R$. 
Then, letting $J_1 = (a_1, \cdots, a_{p-1}, b_{p_1}, a_{p+1},\cdots ,a_s)$ 
and $J_2 = (a_1, \cdots, a_{p-1}, b_{p_2}, a_{p+1},\cdots ,a_s)$ we see by 
induction 
on $s$, 
\[
\del_I^J = \alpha_1\del_I^{J_1} + \alpha_2 \del_I^{J_2}.
\]
Thus,
$D_{\K}(R)$  is generated as a left $D^0_{\K}(R)$-module by the set 
\[
 \{ \del_I^J \mid I =(i_1,\cdots ,i_s) ,J=(a_1,\cdots ,a_s),
i_j \in \N, \textit{ words }a_j \in R, 1\leq i_j \leq n \} \cup 
\{ \lambda_1 \}.
\]
Moreover,
since for 
$a\in R$, 
\[
\lambda_a - \rho_a = \sum_i \del_i^{[a,x_i]}, \textit{ we have }
\rho_a = \lambda_a - \sum_i \del_i^{[a,x_i]}.
\]
From Theorem \ref{theorem-product of dels} and part (2) of Proposition 
\ref{commutator-properties}, we see that
$D_{\K}(R)$ as a left $R$-module (via $\lambda_R$) is generated by the set 
$
 \{ \del_I^J \mid I =(i_1,\cdots ,i_s) ,J=(a_1,\cdots ,a_s),
i_j \in \N, \textit{ words }a_j \in R, 1\leq i_j \leq n \} \cup 
\{ \lambda_1 \}$.

In particular, $D_{\K}^1(R)$ is a left $R$-module generated by
$\{ \del_i ^a \mid a \textit{ is a word in} R \} \cup \{ \lambda_1 \}$.
We suspect that the following proposition is well-known. But we have
not found it published anywhere. 

\begin{proposition}\label{prop-der-free}
The left $R$-module $D^1_{\K}(R)$ is free with basis 
\[
\{ \del_i ^a \mid a \textit{ is a word in } R \} \cup \{ \lambda_1 \}.
\]
\end{proposition}
\begin{proof}
Suppose $\varphi = a\lambda_1+ \sum_{i=1}^n \sum_{j=1}^{r_i} h_{ij}
\del_i^{f_{ij}}   
=0$
for some $a, h_{ij}\in R$, and $f_{ij}$ words in $R$.
Since $\varphi (1) =0$, we have, $a=0$. 
Hence, $\varphi = \sum_{i=1}^n \sum_{j=1}^{r_i} h_{ij}\del_i^{f_{ij}}   =0$.
Now, for each $i\leq n$ and any word $w\in R$, we have
\[
\varphi (wx_i) = \varphi (w)x_i + \sum_{j=1}^{r_i} h_{ij}wf_{ij} =0.
\]
Since $\varphi (w)=0$, we have 
$ \sum_{j=1}^{r_i} h_{ij}wf_{ij} =0$ for every word $w\in R$. 
That is, $\sum_{j=1}^{r_i} h_{ij} \otimes f_{ij}^o \in 
\textit{ Ker }:R\otimes R^o 
\to D_{\K}^0(R)$
of proposition \ref{proposition-0-order}.
That is, $\sum_{j=1}^{r_i} h_{ij} \otimes f_{ij}^o =0$, hence the result.
\end{proof}
\begin{remark}
Note that this result does not generalize to $D^i_{\K}(R)$ for $i\geq 2$. 
For example,
\[
\del_{(1,2)}^{(x_2x_1,1)} - \del_{(1,2)}^{(x_1x_2, 1)} - x_2\del_2^1 +
\del_2^{x_2} =0.
\]
But we can generalize Proposition \ref{prop-der-free} in the following sense.
\end{remark}
\begin{proposition}\label{I-fixedlength-free}
For a fixed $s\geq 1$, the set 
\[
\{ \del_I^J \mid I=(i_1,\cdots ,i_s), J=(a_1,\cdots ,a_s), i_j\leq n, 
\textit{ words }a_j \in R \}
\]
generates a free $D^0_{\K}(R)$-submodule of $D^s_{\K}(R)$. 
\end{proposition}
\begin{proof}
The argument is essentially the same as in the proof of 
Proposition \ref{prop-der-free} and we use induction on $s$.
Suppose $\varphi = \sum_{|I|=s} \sum_{k} \rho_{q_{IJ}^k} 
\lambda_{p_{IJ}^k} \del_I^J =0$ for some
finitely many monomials $q_{IJ}^k, p_{IJ}^k \in R$ with every entry 
of $J$ a word. For any $t\leq n$, and word
$w \in R$, we have
$\varphi( wx_t) = \varphi (w)x_t + \sum_{\{ I \mid i_s=t \}} \sum_k  
\rho_{a_s} \rho_{q_{IJ}^k} \lambda_{p_{IJ}^k} 
\del_{\widehat{I}}^{\widehat{J}} (w)  $ where $\widehat{I} 
= (i_1,\cdots ,i_{s-1})$ and
$\widehat{J} = (a_1,\cdots ,a_{s-1})$, and we use the fact that
$\del_I^J \rho_{x_t} = \rho_{x_t} \del_I^J+ 
\rho_{\del_{i_s}^{a_s}(x_t)}\del_{\widehat{I}} ^{\widehat{J}}$.
Hence, $\varphi =0$ implies that $ \sum_{\{ I \mid i_s=t \}} 
\sum_k  \rho_{a_s} \rho_{q_{IJ}^k} \lambda_{p_{IJ}^k} 
\del_{\widehat{I}}^{\widehat{J}} =0$ and we now appeal to induction to 
complete the proof.
\end{proof}
\begin{proposition}\label{D(R)-span}
 $D_{\K}(R)$ is spanned as a $\K$-vector space by the set
\begin{align*}
 &\{ \del_I^J \mid I =(i_1,\cdots ,i_s) ,J=(a_1,\cdots ,a_s),
i_j \in \N, \textit{ words }a_j \in R, 1\leq i_j \leq n \}\\ 
 &\cup \{ \lambda_a 
\mid a 
\textit{ is a word in } R \}.
\end{align*}
\end{proposition}
\begin{proof}
We already know that
$D_{\K}(R)$ as a left $R$-module (via $\lambda_R$) is generated by the set 
$
 \{ \del_I^J \mid I =(i_1,\cdots ,i_s) ,J=(a_1,\cdots ,a_s),
i_j \in \N, \textit{ words }a_j \in R, 1\leq i_j \leq n \} \cup 
\{ \lambda_1 \}$.

Now, let $a\in R$ be a word, $I=(i_1,\cdots ,i_s)$, $J=(a_1,\cdots ,a_s)$ 
for some 
$s\geq 1$, $i_j \leq n$, and words $a_j \in R$. Let $\widehat{I} = 
(i_2,\cdots ,i_s)$ and
$\widehat{J} = (a_2,\cdots ,a_s)$. Then,
\begin{align*}
[a\del_I^J , x_{i_1}] &= [a,x_{i_1}] \del_I^J + 
aa_{1}\del_{\widehat{I}}^{\widehat{J}}, \textit{ and }\\
[a\del_I^J, x_p] &= [a,x_p]\del_I^J \textit{ for }p\neq i_1.
\end{align*}
Further, $a\del_I^J (1)=0$. That is, letting $(p,I)= (p, i_1,\cdots ,i_s)$ and 
$([a,x_p],J) = ([a,x_p], a_1,\cdots ,a_p)$ for $1\leq p \leq n$, and
$a\cdot J = (aa_1,a_2,\cdots ,a_s)$ we have
\[
a\del_I^J =  \del_I^{a\cdot J}+\sum_{p=1}^n \del_{(p,I)}^{([a,x_p],J_p)}.
\]
\end{proof}
\begin{proposition}\label{goingup}
Let $I_1 = (i_1,\cdots ,i_r), J=(j_1,\cdots ,j_s)$, $A=(a_1,\cdots ,a_r)$,\\
$B=(b_1,\cdots ,b_s)$ with $i_l,j_l \leq n$, 
and $a_l,b_l \in R$. Let $A\cdot w = (a_1,\cdots ,a_rw)$ and 
$w\cdot B = (wb_1,\cdots, b_s)$ for $w\in R$.
Then
\[
\del_{(I ,J)}^{(A\cdot w, B)} - \del_{(I,J)}^{(A, w\cdot B)} = 
\sum_{k=1}^n \del_{(I, k, J)}^{(A, [w,x_k], B)}.
\]
\end{proposition}
\begin{proof}
This is proved by induction on $r$. Note
\begin{align*}
[\del_{(i,J)}^{(a_1\cdot w, B)} , x_i] &= a_1\cdot w \del_J^B 
					= a_1\left( \del_J^{w\cdot B} 
                                          + \sum_{k=1}^n \del_{(k, J)} 
                                          ^{([w,x_k], B)}
						\right)\\
			&\hspace{0.2in} 
                        \text{(proof of Proposition \ref{D(R)-span})};\\
[\del_{(i,J)}^{(a_1\cdot w, B)}, x_t] &= 0 \quad \text{ for }t\neq i, \quad 
\text{and}\quad 
\del_{(i,J)}^{(a_1\cdot w, B)} (1) =0.\\
\text{Hence,} \quad
\del_{(i,J)}^{(a_1\cdot w, B)} &= \del_{(i,J)}^{(a_1,w\cdot B)} + \sum_{k=1}^n 
\del_{(i,k,J)}^{(a_1, [w,x_k], B)}  .
\end{align*}
Now, for $r>1$ and using induction,  
\begin{align*}
[\del_{(I,  J)}^{(A\cdot w,B)}- \del_{(I,J)}^{(A, w\cdot B)}   , x_t] &= 
[ \sum_{k=1}^n \del_{(I, k, J)}^{(A, [a,x_k], B)} , x_t] \quad 
\text{for }t\leq n, 
\text{ and}\\
 \del_{(I,  J)}^{(A\cdot w,B)}- \del_{(I,J)}^{(A, w\cdot B)} (1)&=0.
\end{align*}
Hence the result. 
\end{proof}
\begin{remark}\label{series operators}
For any indexing tuple $I = (i_1,\cdots ,i_r)$, and a word $w\in R$, denote
by $w_I$, the tuple $(1,\cdots ,1,w)$. 
The above proposition implies that every $\varphi \in D_{\K}(R)$ can be 
written in 
the form of a
series 
\[
\varphi = \lambda_a + \sum_{I,w} \alpha_{I,w} \del_I^{w_I}
\quad \text{for $\alpha_{I,w} \in \K$, $a \in R$, and words $w \in R$.}
\]
But we have the following finite and unique description for any $\varphi 
\in D_{\K}(R)$ noting that
$D_{\K}(R)$ as a vector space has a spanning set as given by Proposition 
\ref{D(R)-span}.
\end{remark}
\begin{corollary}\label{D(R)-somebasis}
Let $r\geq 0$ be the minimum natural number for which $\varphi $ is in the 
span of 
the set
$\{ \lambda _a \mid a \text{ is a word in } R \} \cup 
\{ \del_I^A \mid |I| \leq r \}$. Then
for finitely many words $w\in R$ and finitely many $A=(a_1,\cdots ,a_r)$,
	$a,a_l \in R$, 
 $\varphi$ can
be uniquely written
\[
\varphi = \lambda_a + \sum_{|I|<r, w} \alpha_{I,w} \del_{I}^{w_I} + 
\sum_{|I|=r, A} 
\alpha_{I,A}
	\del_I^A. 
\]
\end{corollary}
\begin{proof}
By Propositions \ref{D(R)-span} and \ref{goingup} we see that $\varphi$ can 
indeed
by written in the given fashion. For uniqueness, suppose
\[
\varphi = \lambda_a + \sum_{|I|<r, w} \alpha_{I,w} \del_{I}^{w_I} + 
\sum_{|I|=r, A} 
\alpha_{I,A}
	\del_I^A =0
\]
for finitely many words $w\in R$,  and finitely many $A=(a_1,\cdots ,a_r)$,
with $a,a_l \in R$.
Then, $\varphi (1)=0$ implies $a=0$. Since, for each $I=(i_1,\cdots ,i_t)$ 
with $t<r$,
$\del_I^{w_I}(x_{i_1}\cdots x_{i_t}) = w$, we have $\varphi(x_{i_1}\cdots 
x_{i_t}) = 
\alpha_{I,w}=0$.
Thus, \\
$\varphi = \sum_{|I|=r, A} \alpha_{I,A}
	\del_I^A =0$. Now we refer to Proposition \ref{D(R)-span} to claim 
        that 
        $\alpha_{I,A} = 0$. 
\end{proof}
\begin{remark}\label{Remark:canonical-series}
The proof of the above Collary also proves that the series description as in 
the 
Remark \ref{series operators} is also unique.
We say that $\varphi$ is written in the {\bf finite canonical form} 
(respectively, 
the {\bf power series canonical form}) 
if it is written as described in the above corollary (respectively, remark 
\ref{series operators}).
\end{remark}
\begin{theorem}\label{D(R) is simple}
The algebra $D_{\K}(R)$ is simple. 
\end{theorem}
\begin{proof}
Suppose $I$ is a nonzero ideal of $D(R)$. Let $\varphi \in I$ be written in 
the 
finite canonical form.
Note that $[\varphi , \del_i^1]$ is also in the finite canonical form with 
reduced 
degree.
Thus, taking successive commutators with appropriate $\del_i^1$ we can 
assume that 
$\varphi = \lambda_1 + \sum_{|I|\leq r}\alpha_I \del_I^{\mathbf{1}} $ where 
$\mathbf{1}=(1,1,\cdots ,1)$, unless, $[\varphi , \del_i^1]=0$ $\forall i$. 
In the former case, we can take appropriate commutators with $x_i$ to
arrive at $1\in I$. In the latter case, we refer to part (2) of the 
Proposition 
\ref{commutator-properties}, and consider
$[\varphi, \del_{(i,j)}^{(1,1)}] \in I$ $\forall i,j$. Taking commutator with 
$\del_{(i,j)}^{(1,1)}$ results in an operator of
reduced degree written in the finite canonical form, unless the commutator 
is 0 
$\forall i,j$. In the latter case, we
take commutators with appropriate $\del_{(i,j,k)}^{(1,1,1)}$. Continuing 
thus we 
arrive at the result. 
\end{proof}
\begin{remark}
We believe that the canonical forms can be used to prove that $D_{\K}(R)$ is
a domain when characteristic of $\K$ is zero. We are unable to prove the
same.

We also believe that 
the algebra generated by $D_{\K}^k(R)$ is not all of $D_{\K}(R)$. We have 
seen in 
Remark \ref{remark-Delta-not-equal-to-D} that $\Delta (R) \neq D_{\K}(R)$.
We are not able to prove that, $D_{\K}(R)$ is not finitely generated.
Conjecturally, even when the characteristic of $\K$ is zero,
$\langle D^k_{\K}(R) \rangle \neq D_{\K}(R)$ for any $k\geq 0$. 
\end{remark}
For $A = \K [y_1,\cdots ,y_n]$ the polynomial algebra over a field 
$\K$ of characterisitc $0$ in $n$
variables, and $I$ be an ideal of $A$, 
we have $\mathcal{S}_I / \mathcal{J}_I\cong D_{\K} (A/I)$ as filtered algebras.
Further, $\mathcal{J}_I = ID_{\K}(A)$
and $\mathcal{S}_I$ is the ring generated by $\mathcal{J}_I$ (see chapter 15 
and the references given in section 15.6 of \cite{MR}).
We would like to see an analogue of these 
statements for the free algebra $R = \K \langle x_1,\cdots ,x_n \rangle$, 
$n>1$,
with a two sided ideal $I$.
\begin{proposition}
The natural map 
$\mathcal{S}_I / \mathcal{J}_I \to D_{\K} (R/I)$
of filtered algebras is surjective.
\end{proposition}
\begin{proof}
Let $\overline{a} \in R/I$ denote the image of $a\in R$ and
$\overline{\varphi } \in D_{\K}(R/I)$ denote the image of 
$\varphi \in \mathcal{S}_I$. Note 
$\lambda_{\overline{a}} = \overline{\lambda_a}$
and $\rho_{\overline{a}} = \overline{\rho_a}$.
Hence $D^0_{\K}(R/I)$ is in the image of $\mathcal{S}_I/\mathcal{J}_I$.

Suppose $\eta \in D^t_{\K}(R/I)$ is such that 
$[\eta, \overline{x_i}]= \overline{\psi_i} \in D^{t-1}_{\K}(R/I)$ for 
$t\geq 1$, and $i\leq n$ and $\eta (\overline{1}) =\overline{0}$. Let 
$\varphi \in D^t_{\K}(R)$ be such that 
$[\varphi ,x_i] = \psi_i \in D^{t-1}_{\K}(R)$ and $\varphi (1) =0$. 
Then $\overline{\varphi} = \eta$ (and therefore $\varphi \in \mathcal{S}_I$).
\end{proof}

\section{$\beta$-differential operators on the 
$\Z^n$-graded free algebra.}
\label{S:betafree}
\newcommand{\Dbeta}{D_{\beta}}
\newcommand{\delbeta}[1]{\del_{\beta,#1}}

Let $R = \K \langle x_1,\cdots, x_n \rangle$ 
be the free $\K$-algebra
generated by variables $x_1,\cdots, x_n$.  
The algebra $R$ is $\Z^n$ graded by setting 
degree of $x_i$ to be $e_i = (0,\cdots , 0,1,0,\cdots ,0)$
where $1$ appears in the $i$-th place. For $i,j\leq n$, let
$q_{ij} \in \K ^*$.
We define $\beta :\Z^n \times \Z^n \to \K^*$ by setting
$\beta (e_i, e_j) = q_{ij}$. In case $q_{ij}q_{ji}=1$ for $i\neq j$, we
have $\beta (a,b)\beta (b,a)=1$ for $a\neq b$ which implies that if
$\varphi _1, \varphi_2$ are two left $\beta$-derivations, then so is
$[\varphi_1, \varphi_2]_{\beta}$.

For each homogeneous
$a\in R$, 
and each $i\leq n$, there is a left $\beta$-derivation
$\delbeta{i}^a$
such that
$\delbeta{i}^a (x_j) = \delta_{i,j}a$. Note that  
$d_{\delbeta{i}^a} = d_a - e_i$ 
(recall the notation, $d_m = $ degree of $m$).
Moreover, for $a\in R$, we have the $\beta$-inner derivation,
\[
\lambda_a - \rho^{\beta}_a = \sum_i \delbeta{i}^{[a,x_i]_{\beta}}.
\]
\begin{remark}
In the case when $q_{ij}q_{ji}=q_{ii}=1$, the vector space of
$\beta$-derivations on $R$, denoted $Der_{\beta}(R)$, is a coloured
Lie algebra, which is not simple. There is a surjection from
$Der_{\beta}(R)$ to $Der_{\beta}(\overline{R})$ where
$\overline{R}$ is the quotient algebra of $R$ subject to the relations
$x_ix_j = q_{ij}x_jx_i$ for $i,j\leq n$. The $\beta$-inner derivations
are in the kernel of this surjection. 
\end{remark}
\begin{definition}
For $r=1, I=(i_1), i_1\leq t$, and $A=(a_1)$, for homogeneous $a_1 \in R$
set $\delbeta{I}^{A}
	= \delbeta{i_1}^{a_1}$.
For an $r\in \N, r\geq 2$, 
let $I =(i_1, i_2,\cdots ,i_r)$ be a sequence of natural numbers
$i_j \leq t$ and $A=(a_1,\cdots ,a_r)$ be a sequence of homogeneous
elements from $R$.
Further, let $\widehat{I} = (i_2,\cdots ,i_r)$ and
$\widehat{A} = (a_2,\cdots ,a_r)$. 
Denote by $\delbeta{I}^A \in D_{\beta}^r(R)$ 
the operator which satisfies the 
commutator rules 
\[
[\delbeta{I}^A , x_{i_1}]_{\beta} = 
	a_1\delbeta{\widehat{I}}^{\widehat{A}}, \quad
[\delbeta{I}^A , x_{l}]_{\beta} = 0 \textit{ for }l \neq i_1, \textit{ and }
\quad \delbeta{I}^A(1) =0.
\]
Let $\Delta_{\beta} (R)$ denote the associative subalgebra of $D_{\beta}(R)$
generated by $D^1_{\beta }(R)$ and \\
$\Delta _{\beta}^r(R) = \Delta_{\beta} (R) \cap D_{\beta}^r(R)$ for $r\geq 1$.
\end{definition}
The proofs of the items in the following theorem are similar to those in the 
section \ref{S:free}.
\begin{theorem}\label{gen-D(sfree)}
\begin{enumerate}
\item	The associative, $\Z^n$-graded 
	algebras $R\otimes^{\beta} R^{\beta o}$ and $D^0_{\beta}(R)$ are
	isomorphic. 
	\ignore{Recall from section \ref{S:Prelim} 
	that the multiplications in $R^{\beta ,o}$ and
	$R\otimes R^{\beta ,o}$ respect the bicharacter $\beta$.}
\item 	The $\beta$-centre of $D_{\beta}(R)$ is $\K$.
\item 	A $\beta$-derivation of $R$ which is in $D^0_{\beta}(R)$ is a sum of 
inner
	$\beta$-derivations.

\item	For homogeneous $a, b \in R$ and $i,j$ between $1$ and $n$, we have
	\[
		[\delbeta{i}^a, \lambda_b]_{\beta}= \lambda_{\delbeta{i}^a(b)}
	\]
	Further, when $q_{ij}q_{ji} = 1$ and $q_{ii} = 1$ we have,
	\[
		[\delbeta{i}^a, \rho^{\beta}_b]_{\beta}
		= \rho_{\delbeta{i}^a(b)}^{\beta}
	\]
	and
	\[
		 [\delbeta{i}^a ,\delbeta{j} ^b]_{\beta}
		= \delbeta{j}^{\delbeta{i}^a (b)}
			- \beta (d_a-e_i, d_b-e_j)
			\delbeta{i}^{\delbeta{j}^b (a)}.
	\]

\item
   	Let $r \geq 1$ and let $I$ and $A$ be sequences $I=(i_1,\ldots,i_r)$ 
        and 
        $A=(a_1,\ldots,a_r)$
	 with $i_t$ between $1$ and $n$ and $a_t \in R$. 
	\[
	\delbeta{I}^{A}(x_{t_1}\cdots x_{t_r})
	= 	\begin{cases}
		a_1\cdots a_r &\textit{ if } 
			(t_1, \cdots ,t_r) = I\\ 
		0 &\textit{ if } (t_1, \cdots ,t_r) \neq I 
         	\end{cases}
	\]
	Further, 
	$\delbeta{I}^A(x_{t_1}\cdots x_{t_k}) = 0$ for $k < r$.
\item   
	$[\delbeta{(i_1,\cdots ,i_r)}^{(a_1,\cdots ,a_r)} , \lambda_a]_{\beta} 
	= \sum_{j=2}^r \alpha_j
		\lambda_{\delbeta{(i_1,\cdots ,i_{j-1})}^
                  {(a_1,\cdots ,a_{j-1})} (a)}
			\delbeta{(i_j,\cdots ,i_r)}^{(a_j,\cdots ,a_r)}$ 
	for $r\geq 2$, where $\alpha_j =
        \beta (d_{a_j\cdots a_r} - 
        e_{i_j}-\cdots -e_{i_r}, d_a - 
        e_{i_1}-\cdots -e_{i_j})$. 

	This generalizes the case of 
	$r=1$ from above:  
	$[\del_{\beta, i}^a , \lambda_b]_{\beta} = 
        \lambda_{\del_{\beta, i}^a(b)}$.
	In particular, 
	$[\delbeta{(i_1,\cdots ,i_r)}^{(a_1,\cdots ,a_r)} , \lambda_a]_{\beta} 
	\in D^{r-1}_{\beta}(R)$.
\item  
	$[\delbeta{(i_1,\cdots ,i_r)}^{(a_1,\cdots, a_r)} , 
        \rho_a^{\beta}]_{\beta} 
	= \sum_{j=2}^r \alpha_j
		\rho^{\beta}_{\delbeta{(i_j,\cdots ,i_r)}^{(a_j,\cdots ,a_r)}
                  (a)} 
		  \delbeta{(i_1,\cdots ,i_{j-1})}^{(a_1,\cdots ,a_{j-1})} $    
	when $r\geq 2$, \\
        $q_{ij}q_{ji} = 1$, $q_{ii} = 1$
        and\\
        $
        \alpha_j= \beta (d_{a_1\cdots a_{j-1}} - e_{i_1}-
        \cdots -e_{i_{j-1}}, 
        d_b)
		 \beta (d_{a_1\cdots a_{j-1}}, d_{a_j\cdots a_r}-e_j-
                 \cdots -e_r)
        $.

        This generalizes the 
        case of 
	$r=1$ from above:  
	$[\delbeta{i}^a , \rho^{\beta}_b]_{\beta} = 
        \rho^{\beta}_{\delbeta{i}^a(b)}$.
	In particular, 
	$[\delbeta{(i_1,\cdots ,i_r)}^{(a_1,\cdots ,a_r)} , 
        \rho^{\beta}_a]_{\beta} 
	\in D^{r-1}_{\beta}(R)$.

\item 
	\[
	  \delbeta{I}^A(ab) = \sum_{(J,K) = I} 
          \alpha_{J,K}\delbeta{J}^{A_J}(a) 
          \delbeta{K}^{A_K}(b)
	\]
	where $\alpha_{J,K} = \beta(\degr{\delbeta{K}^A},
        \degr{\delbeta{J}^{A_J}})$. 
\item  
	Let $J=(j_1, \ldots, j_s)$ and $B=(b_1,\ldots,b_s)$.  
	Then there are scalars $\alpha_{\tau}$ such that
     	\[
	  \delbeta{I}^A \delbeta{J}^B = \sum_{\tau \in T_{(r,s)}} 
 	  \alpha_{\tau} \delbeta{\tau((I,J))}^{\tau((A,B))}
 	  + \text{terms of lower order}
	 \]
         Recall that $T_{(r,s)}$ consists of $\tau \in S_{r+s}$ such that
         $\tau$ is increasing on $1,\cdots ,r$ and on $r+1,\cdots ,r+s$.
\item When $q_{ij}q_{ji} = 1$ and $q_{ii} = 1$ for all $i$ and $j$, then
	for $\delbeta{I}^A \in \Dbeta^r$ and $\delbeta{J}^B \in \Dbeta^s$,
	 we have $[\delbeta{I}^A,\delbeta{J}^B]_{\beta} \in \Dbeta^{r+s-1}$.
\item  We have the formula
	\[
        \delbeta{(I,J)}^{(A\cdot w,B)} - \delbeta{(I,J)}^{(A, w \cdot B)}
	  = \sum_{k = 1}^n \alpha_{k}\delbeta{(I,k,J)}^
          {(A,[w,x_k]_{\beta},B)}
          \]
          where $\alpha_k = \beta(\degr{\delbeta{J}^B}, \degr{x_k})$.

\item	For $r \geq 1$,
	Then there are scalars $\alpha_{\sigma}$ such that
	 \[
		\sum_{\sigma \in S_r} \alpha_{\sigma} \delbeta{\sigma(I)}^
                {\sigma(A)} 
		\in \Delta_{\beta}(R)
	\]
	 where $S_r$ is the permutation group on $r$ symbols.

\item
	$\Delta_{\beta}(R) \neq \Dbeta(R)$

\item
	When $q_{ij}q_{ji} = 1$ and $q_{ii} = 1$, the algebra $\Delta_{\beta}
        (R)$ is 
        not simple.

\ignore{\item 
	$\Dbeta^1(R)$ is free as a left $R$-module with basis 
	$\{ \delbeta{i}^a \mid \text{$a$ is a word in $R$}\} \cup 
        \{ \lambda_1 \}$.}
\item
	As a left $\Dbeta^0(R)$-module, $\Dbeta^r(R)$ is generated by the set 
	\[
        \{ \delbeta{I}^A \mid I = (i_1,\ldots,i_r), A=(a_1,\ldots,a_r)\} \cup 
        \{ \lambda_1 \}.
        \]
\item
	The set $\{ \delbeta{I}^A \mid I = (i_1,\ldots,i_r), 
        A=(a_1,\ldots,a_r) 
	\text{ $a_i$ is a word in $R$}\}$ generates a 
	free  $\Dbeta^{0}(R)$-submodule of $\Dbeta^{r}(R)$.
\item
	As a $K$-vector space, $\Dbeta(R)$ is spanned by the set 
   	\begin{align*}
		&
                \{ \delbeta{I}^A \mid I = (i_1,\ldots,i_r), A=(a_1,\ldots,a_r) 
		\text{where each $a_i$ is a word in $R$}\} \\
                &\cup 
		\{ \lambda_a \mid \text{$a$ is a word in $R$} \}
	\end{align*}
\item
	The algebra $\Dbeta(R)$ is simple when $q_{ij} = q_{ji}$ and 
        $q_{ii}=1$.
\item 
  The map described in Proposition \ref{super-S/RtoD} is surjective
  when $R$ is free.
\end{enumerate}
\end{theorem}

\section{Quantum differential operators on the 
  $\Z^n$-graded free algebra.}
\label{S:qfree}
As in the previous section, let $R = \K <x_1,\cdots ,x_n>$ 
be the free $\K$-algebra
generated by variables $x_1,\cdots ,x_n$.  
The algebra $R$ is $\Z^n$ graded by setting 
degree of $x_i$ to be $e_i $, and we 
define $\beta :\Z^n \times \Z^n \to \K^*$ by setting
$\beta (e_i, e_j) = q_{ij}$ for fixed
$q_{ij}\in \K ^*$.

For each $\gamma \in \Gamma$, recall $\sigma_{\gamma} \in D^0_q (R)$ given by
$\sigma_{\gamma}(r) = \beta (\gamma ,d_r) r$ for homogeneous $r$ and extended 
linearly, the grading map.
Let $\Lambda \subset \Aut_{\K}(R)$ be the subgroup generated by
$\{ \sigma_{\gamma} \mid \gamma \in \Gamma \}$. Then we have a surjection 
$(R\otimes R^o)\# \Lambda \to D_q^0(R)$ (recall from Section \ref{S:Prelim} the
surjection $(R\otimes R^o)\# \Gamma \to D_q^o(R)$).
Note that for any homogeneous $\varphi \in \grHom_{\K}(R,R)$ of degree 
$d_{\varphi}$,
 we have
\[
\sigma_{\gamma} \varphi = \beta (\gamma ,d_{\varphi}) \varphi \sigma_{\gamma}.
\]
Hence, every $\varphi \in D^0_q(R)$ can be written as a finite sum,
        $\varphi = \sum_{\sigma \in \Lambda}\psi_{\sigma} \sigma$ for 
        $\psi_{\sigma}
        \in D^0_{\K}(R)$. 

For homogeneous $\varphi, \psi$, let
$[\varphi, \psi]_{\gamma} = \varphi \psi - \beta (\gamma, d_{\psi})\psi 
\varphi$.
The \emph{quantum centre} of $D_q(R)$ is the subalgebra generated by the set
$\{ \text{homogeneous }\varphi \in D_q(R) \mid \exists \gamma \in \Gamma 
\text{ such that }
       [\varphi, \psi]_{\gamma}=0 \forall \psi \in D_q(R) \}$.
\begin{proposition}
The quantum centre of $D_q(R)$ is $\K$.
\end{proposition}
\begin{proof}
Let homogeneous $\varphi $ be in the quantum centre of $D_q(R)$.
Then there exists a $\gamma \in \Gamma$ such that
$[\varphi, \psi]_{\gamma}=0, \forall \psi \in D_q(R)$. In particular,
$[\varphi, \lambda_r]_{\gamma}= 0, \forall r\in R$. Hence,
$\varphi = \rho_{\varphi (1)}\sigma_{\gamma}$. Hence, for any $r\in R$,
$[\varphi, \rho_r]_{\gamma} = \rho_{[\varphi (1), \sigma_{\gamma}(r)]}
\sigma_{\gamma}$.
Thus, $\varphi (1)$ is in the usual centre of $R$.      
\end{proof}
\begin{definition}
For each $i\leq n, \gamma \in \Gamma$, and $a\in R$, denote by 
$\del_{i,\gamma}^a 
\in D^1_q(R)$
defined by $[\del_{i,\gamma}^a , x_j] = \delta_{i,j} \lambda_a 
\sigma_{\gamma}$ and
$\del_{i,\gamma}^a (1)=0$.

For each natural number $r\geq 1$, let $I=(i_1, i_2,\cdots ,i_r)$,  
$K=(\gamma_1, 
\gamma_2,\cdots ,\gamma_r)$,
and  $A=(a_1,\cdots ,a_r)$ for 
$1\leq i_1,\cdots ,i_r \leq n$, $\gamma_1,\cdots ,\gamma_r \in \Gamma$,
and $a_1,\cdots ,a_r \in R$.

When $r=1$, let $\del_{I,K}^A = \del_{i_1,\gamma_1}^{a_1}$.

For $r\geq 2$, let $\del_{I,K}^A \in D_q^r(R)$ be that operator defined by
\[
[\del_{I,K}^A , x_{j}] = \delta_{i_1,j}a_1\del_{\widehat{I}, 
  \widehat{K}}^{\widehat{A}} \sigma_{\gamma_1},\quad
\del_{I,K}^A (1)=0,
\] 
where
$\widehat{I}= (i_2,\cdots ,i_r), \widehat{K}= (\gamma_2,\cdots ,\gamma_r)$, and
$\widehat{A}= (a_2,\cdots ,a_r)$.
\end{definition}
\begin{proposition}
For each $i\leq n$, $a\in R$, and $\gamma \in \Gamma$, the operator 
$\del_{i,\gamma}^a$ is a right skew $\sigma_{\gamma}$-derivation.
\end{proposition}
\begin{proof}
We need to prove that $[\del_{i,\gamma}^a, r] = 
\lambda_{\del_{i,\gamma}^a (r)}\sigma_{\gamma}$ $\forall r\in R$.
We know that this is true for $r \in \{x_1,\cdots ,x_n, 1 \}$.
Assume that the proposition is true for every word of length less than
$k$. Suppose that $rs$ is a word of length $k$. 
Then, 
\begin{align*}
[\del_{i,\gamma}^a ,rs] &= [\del_{i,\gamma}^a,r]s +r[\del_{i,\gamma}^a, s]
                   =\lambda_{\del_{i,\gamma}^a(r)}\sigma_{\gamma} s+r 
                   \del_{i,\gamma}^a (s)\sigma_{\gamma}\\
                   &=\lambda_{\del_{i,\gamma}^a(r)}\sigma_{\gamma}(s)
                   \sigma_{\gamma}+r \del_{i,\gamma}^a (s)
                   \sigma_{\gamma}
                   = (\del_{i,\gamma}^a(r)\sigma_{\gamma}(s) +
                   r\del_{i,\gamma}^a(s)) 
                   \sigma_{\gamma}
\end{align*}
Now, $\del_{i,\gamma}^a (rs) = \del_{i,\gamma}^a(rs1)=rs\del_{i,\gamma}^a(1) + 
[\del_{i,\gamma}^a,rs](1)
=\del_{i,\gamma}^a(r)\sigma_{\gamma}(s) +r\del_{i,\gamma}^a(s)$. 
\end{proof}
\begin{remark}\label{der-mult}
The above proposition shows that 
\[
\del_{i,\gamma }^a \lambda_r - \lambda_r\del_{i,\gamma}^a = 
\lambda_{\del_{\gamma,i}^a (r)} \sigma_{\gamma} 
\quad \textit{and} \quad
\del_{i,\gamma }^a \rho_s - \rho_{\sigma_{\gamma}(s)}\del_{i,\gamma}^a = 
\rho_{\del_{\gamma,i}^a (s)}
\]
We see that $D_q^1(R)$ is generated as a module over $D_q^0(R)$ by 
the right $\sigma_{\gamma}$-derivations $\del_{i,\gamma}^a$. Recall that
$\varphi$ is a left skew $\sigma$-derivation, if and only if
$\varphi \sigma^{-1}$ is
a right skew $\sigma^{-1}$-derivation.  

Let $\Delta_q(R)$ be the subalgebra of $\grHom_{\K}(R,R)$
generated by $D_q^0(R)$ and the operators $\del_{i,\gamma}^a $ for
$i\leq n, a\in R$ and $\gamma \in \Gamma$.

For any $a\in R$, and $\gamma \in \Gamma$, the operator
$\lambda_a - \rho_a \sigma_{\gamma}$ is a left skew-$\sigma_{\gamma}$ 
derivation.
We call such a left skew-$\sigma_{\gamma}$ derivations, an inner left-
$\sigma_{\gamma}$ 
derivation. Following the same proof as in Proposition 
\ref{inner-derivations}, we see
that any left skew-$\sigma_{\gamma}$ derivation of $R$ which is in $D^0_q(R)$
is a sum of inner left-$\sigma_{\gamma}$ derivations.

Similarly,  $\lambda_a \sigma_{\gamma} - \rho_a$ is a right 
skew-$\sigma_{\gamma}$ derivation. Such a right 
skew-$\sigma_{\gamma}$ derivations is called, an inner 
right-$\sigma_{\gamma}$ derivation. Any 
right skew-$\sigma_{\gamma}$ derivation of $R$ which is in $D^0_q(R)$
is a sum of inner right-$\sigma_{\gamma}$ derivations.
\end{remark}
Here are some more generalizations of the usual and $\beta$ differential
operators. The proofs are similar to the corresponding ones in
the section \ref{S:free}.
\begin{theorem}\label{gen-D(qfree)}
  Let $I=(i_1,\cdots, i_r)$, $A=(a_1,\cdots ,a_r)$, and 
  $K=(\gamma_1,\cdots ,\gamma_r)$ as 
  in the definition above. 
\begin{enumerate}
\item   
  \[
  \del_{I, K}^A (x_{t_1}x_{t_2}\cdots x_{t_r}) =
  \begin{cases}
   \left( \prod_{1\leq j<s\leq r}\beta (\gamma_{j}, e_{i_s}) \right)
   a_1a_2\cdots a_r
    &\textit{if } I= (t_1,t_2,\cdots ,t_r)\\
    0 &\textit{if } I\neq (t_1,t_2,\cdots ,t_r)
  \end{cases}
  \]
  Further, $  \del_{I, K}^A (x_{t_1}x_{t_2}\cdots x_{t_k}) =0$
  for $k<r$.
\item	Let $J=(i_{j_1},\cdots ,i_{j_t})$ be a 
  subsequence of $I$. 
  Let $A_{J} = (a_{j_1},\cdots ,a_{j_t})$ and $K_{J} = 
  (\gamma_{j_1},\cdots ,\gamma_{j_t})$  
  respectively represent the corresponding subsequences of $A$ and $K$.
  Further, let $\sigma_{K_J} = \sigma_{\gamma_{j_1}} +\cdots 
  +\sigma_{\gamma_{j_t}}$.
  Then, for $a,b\in R$,
  \[
  \del_{I,K}^A(a\cdot b) = \sum_{(I_1,I_2)=I} \del_{I_1,K_{I_1}}^{A_{I_1}} 
  (a)\cdot 
  \del_{I_2,K_{I_2}}^{A_{I_2}} (\sigma_{K_{I_1}} (b))
  \]
 where $(I_1,I_2)$ denotes the concatenation of $I_1$ and $I_2$.
Therefore, (recall Remark \ref{der-mult})
\[
\del_{I,K}^A \lambda_a - \lambda_a \del_{I,K}^A = \sum_{(I_1,I_2)=I, 
  I_1,I_2 \neq \{ \}} 
  \lambda_{\del_{I_1,K_{I_1}}^{A_{I_1}} (a)}
  \del_{I_2,K_{I_2}}^{A_{I_2}} \sigma_{K_{I_1}} 
\]
\[
\del_{I,K}^A \rho_b - \rho_{\sigma_K(b)} \del_{I,K}^A =
   \sum_{(I_1,I_2)=I, I_1,I_2 \neq \{ \}} \rho_{\del_{I_2, K_{I_2}}^{A_{I_2}} 
     (\sigma_{K_{I_1}}(b)) }
     \del_{I_1,K_{I_1}}^{A_{I_1}}.
\]
\item Let $S_r$ denote the permutation group on $r$ symbols, and
  for any $\tau \in S_r$ let $\tau (I) = (i_{\tau (1)},
  i_{\tau (2)}, \cdots ,i_{\tau (r)})$, and similarly define
  $\tau (A), \tau(K)$. Then, for each $\tau \in S_r$ 
  there exists $\alpha_{\tau} \in \K$
  such that 
  \[
  \sum_{\tau \in S_r} \alpha_{\tau}\del_{\tau (I),\tau (K)}^{\tau (A)} 
  \in \Delta_q(R).
  \]
  In fact, $ \sum_{\tau \in S_r} 
  \alpha_{\tau}\del_{\tau (I),\tau (K)}^{\tau (A)} -
  \del_{i_1,\gamma_1}^{a_1}\del_{i_2,\gamma_2}^{a_2}
  \cdots \del_{i_r,\gamma_r}^{a_r} \in \Delta_q(R) \cap D^{r-1}_q(R)$.
  Moreover, if $d_m$ denotes the degree of the operator
  $\del_{i_m,\gamma_m}^{a_m}$ for $1\leq m \leq r$, then
  the scalar 
  \[
  \alpha_{\tau} = \prod_{1\leq \tau(m)<\tau(n)\leq r} \beta 
  (\gamma_{\tau(m)} ,d_{\tau(n)}). 
  \]
\item  
	Let $J=(j_1, \ldots, j_s)$, $B=(b_1,\ldots,b_s)$, and
        $L=(\delta_1,\cdots ,\delta_s)$.  
	Then there are scalars $\alpha_{\tau}$ such that
     	\[
	  \del_{I,K}^A \del_{J,L}^B = \sum_{\tau \in T_{(r,s)}}
          \alpha_{\tau}\del_{\tau (I,J) , \tau (K,L)}^{\tau (A,B)}
 	  + \text{terms of lower order}
	 \]
         Recall that $T_{(r,s)}$ consists of $\tau \in S_{r+s}$ such that
         $\tau$ is increasing on $1,\cdots ,r$ and on $r+1,\cdots ,r+s$.
         For $1\leq m \leq r+s$, let 
         $d_m =$ degree of $\del_{i_m, \gamma_m}^{a_m}$ and
         $\eta_m =  \gamma_m$ for $m\leq r$
         and $d_m =$ degree of 
         $\del_{j_{m-r}, \delta_{m-r}}^{b_{m-r}}$ for $m> r$ and
         $\eta_m = \delta_{m-r}$.
         Then 
         \[
         \alpha_{\tau} = \prod_{1\leq \tau(m)<\tau(n)\leq r+s} \beta 
         (\eta_{\tau(m)} ,d_{\tau(n)}). 
         \]
\item   We have the following formula
          \[
          a\del_{I,K}^A - \del_{I,K}^{a\cdot A} =
          \sum_{k=1}^n \del_{(k,I), (0,K)}^{([a,x_i],A)};
          \]
          here, $(a\cdot A) = (aa_1,a_2,\cdots ,a_r)$, 
          $(k,I)=(k,i_1,\cdots ,i_r)$, \\
          $(0,K) = (0,\gamma_1,\cdots ,\gamma_r)$,
          and $([a,x_k],A) = ([a,x_k],a_1,\cdots ,a_r)$.

         With the same notations as used above, or in previous sections,
         \[
         \del_{(I,J),(K,L)}^{(A\cdot w, B)}-
                  \del_{(I,J),(K,L)}^{(A,w\cdot B)}
                  = \sum_{k=1}^n \del_{(I,k,J),(K,0,L)}^{(A,[w,x_k],B)}
         \]
\item  As a left $D_q^0(R)$-module, $D_q^r(R)$ is generated by
  the set 
  \[
  \{ \del_{I,K}^A \mid I=(i_1,\cdots ,i_r),K=(\gamma_1,\cdots,
  \gamma_r), A=(a_1,\cdots ,a_r) \} \cup \{ \lambda_1 \}.
  \]
  \item 
  The map described in Proposition \ref{q-S/RtoD} is surjective
  when $R$ is free.
\end{enumerate}
\end{theorem}
\begin{remark}
\begin{enumerate}
\item 
Note that $D_{\K}(R) \subset D_q(R)$ and we have a map
\[
D_{\K}(R) \# \Lambda \to D_q(R).
\] 
In general, this map need not be
surjective. Indeed, in \cite{IM} the generators of the algebra $D_q (\K [x])$
have been described over a field of characteristic $0$.  The polynomial algebra
$\K [x]$ is $\Z$-graded, and $\beta : \Z \times \Z \to \K^*$ is
given by $\beta (n,m)=q^{nm}$ for $q$ a transcendental element over $\Q$, and
$\Q (q) \subset \K$.  Then $D_q( \K [x])$ is generated by the set 
$\{ \lambda_x, \del^{\beta}, \del,
\del^{\beta^{-1}} \}$, where $\del (x^n)=nx^{n-1}$, 
$\del^{\beta}(x^n) = \left( \dfrac{q^n-1}{q-1} \right) x^{n-1}$, 
$\del^{\beta^{-1}} (x^n) = \left( \dfrac{q^{-n}-1}{q^{-1}-1} \right) x^{n-1}$,
 for
$n\geq 1$ and $\del(1) = \del^{\beta}(1)=\del^{\beta^{-1}} =0$. The algebra
$D_{\K}(\K [x])$ is the first Weyl algebra, with generators $\{ \lambda_x, 
\del \}$.
We see that $\del^{\beta} \notin D_{\K}(\K [x])\# \Z$. Indeed, by degree 
considerations, 
if $\del^{\beta} \in  D_{\K}(\K [x])\# \Z$, then 
$\del^{\beta} =  \alpha \del \sigma$ for some $\alpha \in \K$, which is not
possible.
\item Corollary \ref{commuting dels}
  and part (10) of Theorem \ref{gen-D(sfree)} do not generalize to the
  case of quantum differential operators. For example, consider 
  $D_q( \K <x_1,x_2>)$ when $q_{11} = q_{22}=1$ and
  $q_{21} = q_{12}^{-1} = q$ for $q \in \K^*$, where $q$ is transcendental 
  over $\Q$, and $\Q (q)\subset \K$.  
  Note that $\del^1_{1,e_1}, \del^1_{1,e_2} \in D_q^1( \K <x_1,x_2>)$, with
  $\del^1_{1,e_1}(x_1^n) = nx_1^{n-1}$, while
  $\del^1_{1,e_2}(x_1^n) = (1+q+\cdots, q^{n-1})x_1^{n-1}$ for $n\geq 1$.
  But $[\del^1_{1,e_1} , \del^1_{1,e_2}]_{\gamma} \notin D^0_q(\K <x_1,x_2>)$ 
  for any $\gamma \in \Gamma$, which can be seen by degree considerations.
\item The question whether the algebra $D_q(R)$ is simple seems to be a
  difficult question to address. Conjecturally, we believe that $D_q(R)$ is
  simple. The algebra $\Delta_q(R)$ we expect to be not simple, as has
  been already checked when $\beta \equiv 1$.
\end{enumerate}
\end{remark}

\bibliographystyle{amsplain}

\end{document}